\documentclass{article}
\usepackage[margin=2cm]{geometry}
\usepackage[affil-it]{authblk}
\usepackage[utf8]{inputenc}
\usepackage{amsmath,amsthm,amssymb}
\usepackage{nccmath}
\usepackage{float}
\usepackage{color}
\usepackage{tikz}
\usepackage{keyval}
\usepackage{ifthen}
\usepackage{booktabs}
\usepackage{epsf}
\usepackage{psfrag}
\usepackage{wrapfig}
\usepackage{color,rotating,psfrag,amsfonts}
\usepackage{todonotes}
\usepackage[dvips]{epsfig}
\usepackage{subfigure}
\usepackage{adjustbox}

\usepackage{makecell, caption, booktabs, multirow, siunitx}

\newcommand{\E}[1]{\ensuremath{\mathbb{E}\left(#1\right)}}
\newcommand{\var}[1]{\ensuremath{\mathbb{V}\left( #1\right)}}

\newcommand{\co}{$\mathrm{CO}_2$ }

\newcommand{\sdisc}{s_{\text{d}}}

\newcommand{\perm}[1][]{\ensuremath{k_{#1}}}

\newcommand{\upscaledPc}{P_{\text{c}}}

\newcommand{\upscaledsat}[1][]{\ensuremath{S_{#1}}}

\usepackage{tikz-network}
\tikzset{network x offset/.initial=1ex,network y offset/.initial=1ex,
adjust size/.style={minimum width=width("\vertex@Label")+2*\pgfkeysvalueof{/tikz/network x offset},
minimum height=height("\vertex@Label")+2*\pgfkeysvalueof{/tikz/network y offset}}}

\makeatletter
\newif\ifcuboidshade
\newif\ifcuboidemphedge

\tikzset{
  cuboid/.is family,
  cuboid,
  shiftx/.initial=0,
  shifty/.initial=0,
  dimx/.initial=3,
  dimy/.initial=3,
  dimz/.initial=3,
  scale/.initial=1,
  densityx/.initial=1,
  densityy/.initial=1,
  densityz/.initial=1,
  rotation/.initial=0,
  anglex/.initial=0,
  angley/.initial=90,
  anglez/.initial=225,
  scalex/.initial=1,
  scaley/.initial=1,
  scalez/.initial=0.5,
  front/.style={draw=black,fill=white},
  top/.style={draw=black,fill=white},
  right/.style={draw=black,fill=white},
  shade/.is if=cuboidshade,
  shadecolordark/.initial=black,
  shadecolorlight/.initial=white,
  shadeopacity/.initial=0.15,
  shadesamples/.initial=16,
  emphedge/.is if=cuboidemphedge,
  emphstyle/.style={thick},
}

\newcommand{\tikzcuboidkey}[1]{\pgfkeysvalueof{/tikz/cuboid/#1}}

\newcommand{\tikzcuboid}[1]{
    \tikzset{cuboid,#1} 
  \pgfmathsetlengthmacro{\vectorxx}{\tikzcuboidkey{scalex}*cos(\tikzcuboidkey{anglex})*28.452756}
  \pgfmathsetlengthmacro{\vectorxy}{\tikzcuboidkey{scalex}*sin(\tikzcuboidkey{anglex})*28.452756}
  \pgfmathsetlengthmacro{\vectoryx}{\tikzcuboidkey{scaley}*cos(\tikzcuboidkey{angley})*28.452756}
  \pgfmathsetlengthmacro{\vectoryy}{\tikzcuboidkey{scaley}*sin(\tikzcuboidkey{angley})*28.452756}
  \pgfmathsetlengthmacro{\vectorzx}{\tikzcuboidkey{scalez}*cos(\tikzcuboidkey{anglez})*28.452756}
  \pgfmathsetlengthmacro{\vectorzy}{\tikzcuboidkey{scalez}*sin(\tikzcuboidkey{anglez})*28.452756}
  \begin{scope}[xshift=\tikzcuboidkey{shiftx}, yshift=\tikzcuboidkey{shifty}, scale=\tikzcuboidkey{scale}, rotate=\tikzcuboidkey{rotation}, x={(\vectorxx,\vectorxy)}, y={(\vectoryx,\vectoryy)}, z={(\vectorzx,\vectorzy)}]
    \pgfmathsetmacro{\steppingx}{1/\tikzcuboidkey{densityx}}
  \pgfmathsetmacro{\steppingy}{1/\tikzcuboidkey{densityy}}
  \pgfmathsetmacro{\steppingz}{1/\tikzcuboidkey{densityz}}
  \newcommand{\dimx}{\tikzcuboidkey{dimx}}
  \newcommand{\dimy}{\tikzcuboidkey{dimy}}
  \newcommand{\dimz}{\tikzcuboidkey{dimz}}
  \pgfmathsetmacro{\secondx}{2*\steppingx}
  \pgfmathsetmacro{\secondy}{2*\steppingy}
  \pgfmathsetmacro{\secondz}{2*\steppingz}
  \ifthenelse{\equal{\dimx}{1}}
    {\foreach \x in {\steppingx,...,\dimx}}
    {\foreach \x in {\steppingx,\secondx,...,\dimx}}
  {     \ifthenelse{\equal{\dimy}{1}}
    {\foreach \y in {\steppingy,...,\dimy}}
    {\foreach \y in {\steppingy,\secondy,...,\dimy}}
    {   \pgfmathsetmacro{\lowx}{(\x-\steppingx)}
      \pgfmathsetmacro{\lowy}{(\y-\steppingy)}
      \filldraw[cuboid/front] (\lowx,\lowy,\dimz) -- (\lowx,\y,\dimz) -- (\x,\y,\dimz) -- (\x,\lowy,\dimz) -- cycle;
    }
    }
    \ifthenelse{\equal{\dimx}{1}}
    {\foreach \x in {\steppingx,...,\dimx}}
    {\foreach \x in {\steppingx,\secondx,...,\dimx}}
  { \ifthenelse{\equal{\dimz}{1}}
    {\foreach \z in {\steppingz,...,\dimz}}
    {\foreach \z in {\steppingz,\secondz,...,\dimz}}
    {   \pgfmathsetmacro{\lowx}{(\x-\steppingx)}
      \pgfmathsetmacro{\lowz}{(\z-\steppingz)}
      \filldraw[cuboid/top] (\lowx,\dimy,\lowz) -- (\lowx,\dimy,\z) -- (\x,\dimy,\z) -- (\x,\dimy,\lowz) -- cycle;
        }
    }
    \ifthenelse{\equal{\dimy}{1}}
    {\foreach \y in {\steppingy,...,\dimy}}
    {\foreach \y in {\steppingy,\secondy,...,\dimy}}
  { \ifthenelse{\equal{\dimz}{1}}
    {\foreach \z in {\steppingz,...,\dimz}}
    {\foreach \z in {\steppingz,\secondz,...,\dimz}}
    {   \pgfmathsetmacro{\lowy}{(\y-\steppingy)}
      \pgfmathsetmacro{\lowz}{(\z-\steppingz)}
      \filldraw[cuboid/right] (\dimx,\lowy,\lowz) -- (\dimx,\lowy,\z) -- (\dimx,\y,\z) -- (\dimx,\y,\lowz) -- cycle;
    }
  }
  \ifcuboidemphedge
    \draw[cuboid/emphstyle] (0,\dimy,0) -- (\dimx,\dimy,0) -- (\dimx,\dimy,\dimz) -- (0,\dimy,\dimz) -- cycle;%
    \draw[cuboid/emphstyle] (0,\dimy,\dimz) -- (0,0,\dimz) -- (\dimx,0,\dimz) -- (\dimx,\dimy,\dimz);%
    \draw[cuboid/emphstyle] (\dimx,\dimy,0) -- (\dimx,0,0) -- (\dimx,0,\dimz);%
    \fi

    \ifcuboidshade
    \pgfmathsetmacro{\cstepx}{\dimx/\tikzcuboidkey{shadesamples}}
    \pgfmathsetmacro{\cstepy}{\dimy/\tikzcuboidkey{shadesamples}}
    \pgfmathsetmacro{\cstepz}{\dimz/\tikzcuboidkey{shadesamples}}
    \foreach \s in {1,...,\tikzcuboidkey{shadesamples}}
    {   \pgfmathsetmacro{\lows}{\s-1}
        \pgfmathsetmacro{\cpercent}{(\lows)/(\tikzcuboidkey{shadesamples}-1)*100}
        \fill[opacity=\tikzcuboidkey{shadeopacity},color=\tikzcuboidkey{shadecolorlight}!\cpercent!\tikzcuboidkey{shadecolordark}] (0,\s*\cstepy,\dimz) -- (\s*\cstepx,\s*\cstepy,\dimz) -- (\s*\cstepx,0,\dimz) -- (\lows*\cstepx,0,\dimz) -- (\lows*\cstepx,\lows*\cstepy,\dimz) -- (0,\lows*\cstepy,\dimz) -- cycle;
        \fill[opacity=\tikzcuboidkey{shadeopacity},color=\tikzcuboidkey{shadecolorlight}!\cpercent!\tikzcuboidkey{shadecolordark}] (0,\dimy,\s*\cstepz) -- (\s*\cstepx,\dimy,\s*\cstepz) -- (\s*\cstepx,\dimy,0) -- (\lows*\cstepx,\dimy,0) -- (\lows*\cstepx,\dimy,\lows*\cstepz) -- (0,\dimy,\lows*\cstepz) -- cycle;
        \fill[opacity=\tikzcuboidkey{shadeopacity},color=\tikzcuboidkey{shadecolorlight}!\cpercent!\tikzcuboidkey{shadecolordark}] (\dimx,0,\s*\cstepz) -- (\dimx,\s*\cstepy,\s*\cstepz) -- (\dimx,\s*\cstepy,0) -- (\dimx,\lows*\cstepy,0) -- (\dimx,\lows*\cstepy,\lows*\cstepz) -- (\dimx,0,\lows*\cstepz) -- cycle;
    }
    \fi 

  \end{scope}
}

\makeatother

\title{Copula modeling and uncertainty propagation in field-scale simulation of \co fault leakage}

\author[1]{Per Pettersson\thanks{Corresponding author: per.pettersson@norceresearch.no}}
\author[2]{Eirik Keilegavlen}
\author[1]{Tor Harald Sandve}
\author[1,2]{Sarah Gasda}
\author[3]{Sebastian Krumscheid}

\affil[1]{NORCE Norwegian Research Centre, N-5838 Bergen, Norway}
\affil[2]{University of Bergen, N-5020 Bergen, Norway}
\affil[3]{Karlsruhe Institute of Technology (KIT), 76131 Karlsruhe, Germany}

\begin{document}

\maketitle

\begin{abstract}
Subsurface storage of \co is an important means to mitigate climate change, and the North Sea hosts considerable potential storage resources. To investigate the fate of \co over decades or even centuries in vast reservoirs, numerical simulation based on realistic models is essential. However, faults and other complex geological structures introduce modelling challenges as their effects on storage operations are subject to high uncertainty due to limited data. In this work, we present a computational framework for forward propagation of uncertainty, including stochastic upscaling and copula representation of multivariate distributions for a \co storage site model with faults. The Vette fault zone in the Smeaheia formation in the North Sea is used as a test case.
The stochastic upscaling method leads to a reduction of the number of stochastic dimensions and the cost of evaluating the reservoir model. A viable model that represents the upscaled data needs to capture dependencies between variables, and allow sampling. Copulas provide representation of dependent multidimensional random variables and a good fit to data, allow fast sampling, and coupling to the forward propagation method via independent uniform random variables. The non-stationary correlation within some of the upscaled flow function are accurately captured by a data-driven transformation model.
The uncertainty in upscaled flow functions and other uncertain parameters are efficiently propagated to uncertain leakage estimates
using numerical reservoir simulation of a two-phase system including \co and brine.
The expectations of leakage are estimated by an adaptive stratified sampling technique, where samples are allocated sequentially and successively concentrated to regions of the parameter space to greedily maximize variance reduction. We demonstrate cost reduction compared to standard Monte Carlo of one or two orders of magnitude for simpler test cases with only fault and reservoir layer permeabilities assumed uncertain, and factors 2-8 cost reduction for stochastic multi-phase flow properties and more complex stochastic models.
\end{abstract}

\section{Introduction}
\co capture and storage (CCS) is anticipated to play a major role in the target reduction of greenhouse gas emissions described in~\cite{Krevor_etal_23}. By 2050, estimates indicate that CCS will contribute with 13\% of the total cumulative emissions reduction~\cite{Ringrose_Meckel_19}. The Norwegian part of the North Sea has an estimated storage potential of 70 gigatons of carbon dioxide (CO$_2$)~\cite{Halland_etal_11}. To put the storage capacity into perspective, the total European Union's \co emission in 2019 were 3.29 gigatons~\cite{UNFCCC}. Hence, the North Sea reservoirs are anticipated to be an important means to obtain permanent storage of sufficient amounts of CO$_2$. To fully exploit the potential of the North Sea reservoirs, it is essential to perform numerical simulations to assess the viability of prospective storage sites. Complex physics and associated uncertainties should be taken into account, for instance regarding faults and fractures that are present in potential storage sites.

Faults are complex geological structures with heterogeneous properties that are associated with challenges in data acquisition. Indeed, faults can act as sealing barriers~\cite{Pei_etal_15}, but can also be associated with leakage risk~\cite{Wu_etal_21}, and the permeability within a fault can vary by orders of magnitude~\cite{Sperrevik_etal_02}. The structural properties, e.g., fault thickness, are typically highly varying~\cite{Childs_etal_90}.
Petrophysical properties are often not directly measured but estimated via other properties, e.g., permeability as a function of shale gouge ratio~\cite{Yielding_etal_97}.
The local fault properties are known only to a limited extent. Seismic imaging cannot resolve the smallest scales, and measurements typically come with uncertainty. Even if increased seismic resolution were possible, it may still not be possible to represent all small-scale features within the scale of the grid discretization possible to achieve in numerical simulation. An important source of uncertainty that is often not taken into account is the frequency of throw variability along a fault, leading to underestimated uncertainty effects~\cite{Manzocchi_etal_08}.
Fluid flow in the fault interacts with flow in adjacent formations, including both the storage formation and potential other formations that are being used in subsurface operations. 
Similarly, uncertainties in fault and formation properties will interact and must therefore be studied together.
Previous work on the Vette fault zone in the North Sea, to be investigated as a test case in this paper, include the work in~\cite{Michie_etal_21}, where the effects on sealing from seismic interpretation methodology were analyzed.
A probabilistic fault system stability investigation with four clay smearing scenarios was carried out in~\cite{Rahman_etal_21}, where the probability for system failure was estimated to be between $10^{-4}$ and $10^{-3}$. A comparison between a stochastic geocellular model and a fault smear model for the Vette Fault Zone was carried out in~\cite{Bjornaraa_etal_22}, demonstrating limited potential for fluid flow.

Comprehensive uncertainty quantification of fault properties is paramount; deterministic base-case simulations often come with order(s) of magnitude errors, and improved measurements of a single source of uncertainty have a limited effect on the accuracy of predictions~\cite{Freeman_etal_08}. While a stochastic fault model in itself does not reduce the uncertainty, it provides a quantitative measure and provide a means to improve predictions.
Stochastic modeling of faults amounts to characterization and representation of both the fault core and the surrounding damage zone. 
First, a structural model needs to be established, for instance using a fault facies model where the fault is divided into discrete geometric objects, each with its own petrophysical properties and different dominating features~\cite{Braathen_etal_09}. The properties of the damage zone of the fault structure is typically dominated by the spatial distribution of deformation bands, and any model of those should honor empirical data~\cite{Schueller_etal_13, Berge_etal_22}. 
Simulation methods for the spatial distribution of lenses in the fault core and deformation bands in the damage zone  include the work in~\cite{Kolyukhin_Tveranger_15}.
In the current work, we consider uncertainty of fault properties with regard to pressure communication, fluid flow, and leakage. Geomechanical aspects will not be considered, although we recognize that fault mechanics is an equally important area of uncertainty to quantify.


\section{Conceptual Approach}
\begin{figure}[H]
\centering
{\includegraphics[width=0.96\textwidth]{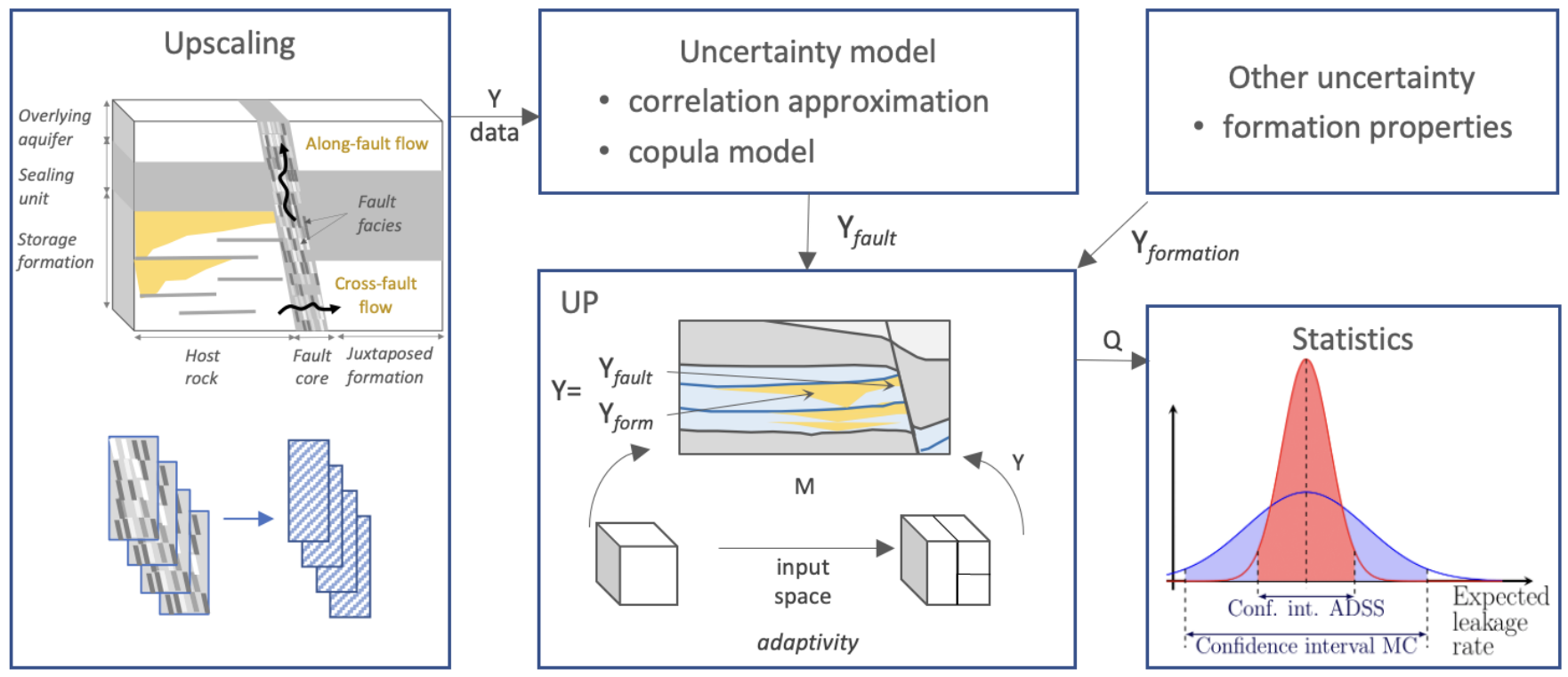}
\label{fig:SGR_data}
}
\caption{Overview of methodological framework. A stochastic geocellular fault model is upscaled (left panel) and a stochastic model is fitted to the upscaled samples (upper center panel). Together with models for other sources of uncertainty, the uncertainties are adaptively propagated through a physical model via random-space stratification (lower middle panel). Expectations of quantities of interest are anticipated to have lower variance compared to standard Monte Carlo estimates (right panel).}
\label{fig:overview_fw}
\end{figure}

In the context of large-scale \co storage simulation, we are thus faced by complex physical problems subject to a large number of possibly dependent uncertainties  interacting with significant physical complexity.
The present work presents a framework for handling this situation, illustrated in Figure~\ref{fig:overview_fw} and explained below. 
To fix notation, 
we aim to perform uncertainty propagation ($U\!P$) for a physical model, denoted by $M$, subject to uncertainties. Specifically, the considered physical model $M$ is a dynamic simulation model that depends on physical parameters $\eta=(\eta_1,...,\eta_d)$, e.g., permeability and relative permeabilities. Furthermore, $M$ will also depend on stochastic variables $Y=(Y_1,..., Y_n)$ to account for the uncertainties affecting the physical problem. These stochastic and physical parameters' dependency is highlighted by the notation $M(Y,\eta)$. 
This distinction will facilitate the narrative and make it clear whether the current focus is on physical or stochastic models.
However, the two are clearly intertwined in that a physical parameter $\eta_i$ can be uncertain (say, a fault permeability is a function of unknown fault properties), which we write as $\eta_i=\eta_i(Y)$.
For a particular realization of $Y$ and $\eta$, the model $M(Y,\eta)$ is used to compute a prediction (e.g., leakage rates) of the physical system. Such an output quantity, denoted $Q$, also depends on $Y$ and $\eta$. Our $U\!P$ approach in this work focuses on computing statistical properties $\mu_Q$, in particular the expected value, of $Q\equiv Q(Y,\eta)$, in short
\begin{equation}
\label{eq:high_level_form}
    \mu_Q = U\!P(M(Y, \eta)).
\end{equation}

While dynamic simulation models for fluid flow through faults in theory can be coupled with reservoir flow simulation, there are
three main reasons to refrain from this path:
First, resolving the complex geometry of the fault will lead to a prohibitively high computational cost of the dynamic model $M$.
Second, instead of randomly sampling existing data as input to the dynamic model $M$, it is often more efficient to first assign or fit a stochastic model for the input distributions.
Third, for scarce data, as will often be the case for geological uncertainty, it may be desirable to add further stochastic modeling assumptions that also need to enter the model, e.g., an assumed but unobserved upper bound on data.
Therefore, in this work we have chosen to treat uncertainty in fault parameters separately from those of the reservoir, 
as is indicated in Figure \ref{fig:overview_fw}, where the two are represented through random variables $Y_{\text{fault}}$ and $Y_{\text{formation}}$, respectively. 

As will be explained in detail in section~\ref{sec:stochastic_faults}, we express fault-related uncertainty through upscaled random variables $Y_{\text{fault}}$ as illustrated in the left panel in Figure \ref{fig:overview_fw}. 
While the upscaling itself follows a standard approach~\cite{rabinovich2016analytical}, which is described in section~\ref{sec:twophase_upscaling}, the coarse-scale flow functions are represented by a stochastic model employing copulas, see section \ref{sec:copula} for details.
Copulas comprise a statistical modeling framework that is suitable for representation and analysis of multivariate (in particular dependent) data. Since the seminal work in~\cite{Joe_96}, building upon the results in~\cite{Sklar_59}, copulas have been used in a wide range of  applications, see~\cite{Bhatti_Do_19} for a list of examples.
There are different parametric families of copulas, each with their own tail dependence, i.e, modeling to what extent an extreme event in one variable is correctly paired with an extreme event in another variable. This is an important feature in the context of \co storage where extreme behavior, e.g., leakage, is of particular interest~\cite{Elsheikh_etal_14, Pettersson_etal_22}.

The copula representation of $Y_{\text{fault}}$ together with a model for $Y_{\text{formation}}$, which we will represent by simple approaches as detailed in section \ref{sec:results}, form the input to the uncertainty propagation (see Figure~\ref{fig:overview_fw}). 
One can distinguish between random sampling methods or statistical methods on the one hand, and deterministic uncertainty propagation methods on the other hand. The former includes variants of Monte Carlo methods, e.g., Latin hypercube sampling~\cite{McKay_etal_79}, multilevel Monte Carlo~\cite{Giles_2015}, and importance sampling~\cite{Kloek_vanDijk_78}, while the latter family of methods includes Quasi-Monte Carlo~\cite{LEcuyer_Lemieux_02} and generalized polynomial chaos via projection~\cite{Xiu_Karniadakis_02} or collocation~\cite{Babuska_etal_07}. 
Whether statistical or deterministic, there is typically a trade-off between robustness and theoretical accuracy.
The computational complexity of the different groups of methods grows more or less with the number of random dimensions and, to some extent, the dependencies between random inputs.  
For the current work, non-intrusive random sampling methods are attractive as they do not require modifications of existing physical solvers, and they typically scale well with the number of random dimensions. In addition, the UP method selected should not impose regularity requirements on the quantity to be approximated.

For uncertainty propagation from random inputs to the quantities of interest, i.e., $\mu_Q$ introduced in Eq.~\eqref{eq:high_level_form}, we will use the adaptive stratified sampling method presented in~\cite{Pettersson_Krumscheid_22},
where samples are adaptively concentrated to regions in input random space that has a large impact on the quantity of interest. The method is sampling-based and does not assume apriori knowledge about what inputs have the greatest impact. It performs best when the variability is localized to a small subset of random variables, or the variance is localized to a subdomain that may not have smaller dimension than the full stochastic space. Related methods that builds upon adaptive sampling and stratification include~\cite{Etore_etal_11,Shields_etal_15,Shields_16, Krumscheid_Pettersson_23}.
The adaptive stratified sampling method employed here is further described in section~\ref{sec:adaptive_sampling}.

For the dynamic model $M$ in Eq.~\eqref{eq:high_level_form}, we utilize existing simulation tools that model large-scale \co storage and thus treat the simulator as a black-box model.
Commercial simulators like Eclipse and Intersect from SLB, and GEM from the CMG group, have dedicated \co modules, as have open-source simulators like DuMuX~\cite{flemisch2011dumux}, GEOX~\cite{gross2021geosx}, MRST~\cite{lie2019introduction}, TOUGH2~\cite{pruess1991tough2} and OPM-Flow~\cite{opmFlowManual}. 
As developing simulation tools is not the focus of this study we use the OPM-Flow simulator, as this has a dedicated \co module, supports standard industry input/output formats, it is fast, and free to use~\cite{sandve2022simulators, sandve2021convective}. 
The latter points are critical for uncertainty quantification, as we will need to run a large number of simulations with varying inputs $\eta$ and $Y$. It is noteworthy that OPM-Flow can account for all relevant fluid flow processes involved in field-scale \co storage in faulted reservoirs, and connected aquifers: \co injection, migration and trapping within the storage formation, along and cross-fault pressure communication and fluid flow.

We test our methodology consisting of stochastic upscaling together with copula representation and further stochastic modeling, and adaptive stratified sampling, on the Smeaheia formation in the North Sea, where the framework will be applied to quantify leakage from the Vette fault zone. 
We base the physical model on the reservoir data-set openly available from CO2Share (https://co2datashare.org/), 
but extend it to include aquifers that are connected to the reservoir through vertical and horizontal flow in the Vette fault zone.
The upscaled stochastic two-phase flow functions specific to the Vette fault are constructed in section \ref{sec:stochastic_faults}. 
These functions depend on $Y_{\text{fault}}$, that is, they reflect uncertainty in the fine-scale distribution of facies within the fault based on available data and a standard fault facies model.
Uncertainty in the formation parameters, $Y_{\text{formation}}$, is introduced in section \ref{sec:results}, where we also present simulations probing the performance of the stochastic framework. The objectives of this work includes demonstration of how the upscaled stochastic flow functions can be included in a reservoir model, and that our sampling strategy gives substantial improvements compared to standard Monte Carlo methods. The methodological components are not limited to the models selected, and hence more widely applicable.
 

\section{Stochastic Fault Model}
\label{sec:stochastic_faults}
Here we describe a stochastic modeling framework for the fault properties on the fine-scale and how these quantities can be upscaled for use as stochastic input for the reservoir simulator. This amounts both to derive parameters for a physically coarser scale, and simultaneously reduce the stochastic complexity by reducing the number of random parameters. While the derivation uses advanced stochastic modeling, the goal is effective stochastic flow functions where each realization is compatible with and has a format similar to what would have been expected in a deterministic model. 

We split the stochastic fault modeling into three stages:
First, the fine-scale geology within the fault is described in terms of a geocellular fault facies model with a stochastic description of the distribution of facies, as well as of their flow properties as we describe in section \ref{sec:stochastic_fault_facies}.
Next, we generate realizations of the fine-scale geometry and compute upscaled permeabilities, and capillary and relative permeability functions valid for the entire fault, see sections \ref{sec:stochastic_fault_facies} and \ref{sec:twophase_upscaling}, respectively.
The final stage entails the representation of the upscaled quantities in a format suited for the stochastic sampling for the field-scale model described in section~\ref{sec:adaptive_sampling}.
In practice, this means the coarse-scale models should
be expressed in terms of independent stochastic variables with known distributions, which is the assumed format for most uncertainty propagation methods.

The coarse-scale absolute permeability can relatively simply be brought into a suitable form, again see section \ref{sec:stochastic_fault_facies}.
For the two-phase flow functions the situation is more challenging, partly reflecting the additional complexity of stochastic functions compared to that of constants (e.g., absolute permeability):
Upscaling from fine-scale realizations produces point clouds of possible two-phase flow properties, which can be represented by equivalent probability distributions that consider only coarse-scale properties, as shown in section~\ref{sec:reduced_stochastic_model}.
This removes the link to the fine-scale model, however, the variables in the stochastic model are statistically dependent. These dependencies can be incorporated into the modelling framework by using so-called copulas, as is done in sections \ref{sec:copula} and \ref{sec:copula_fit}.

While arguably complicated, the presented approach has two main advantages compared to applying a more direct coupling between the stochastic fine-scale geology and the field-scale simulations:
First, properties such as dimension, variability, and heterogeneity of the stochastic fault facies model will not be seen directly by the field-scale reservoir model.
Thus, advanced stochastic descriptions of fault facies can be applied and the added computational effort will only be seen in an offline stage, which is often cheap compared to field-scale simulations.
Second, the modularization allows for applying different level of sophistication for the individual components. For instance, the analytical upscaling described in sections \ref{sec:stochastic_fault_facies} and \ref{sec:twophase_upscaling} could have been replaced by numerical upscaling.
Although this may also necessitate modifying the coarse-scale stochastic modeling to avoid incurring modeling errors in that step (a more complex stochastic model may be needed), there is still further significant flexibility in how our framework is implemented. For instance, a different fine-scale model can be used as long as it is possible to produce independent samples, and the coarse-scale parameters are not assumed to be independent or belong to any specific class of distributions.
Indeed, we expect this flexibility to be particularly useful in the stochastic modeling, since the coarse-scale properties to be represented may vary significantly with the fine-scale model. 

\subsection{Stochastic Modeling of Fault Facies Distribution}
\label{sec:stochastic_fault_facies}
As a simple non-trivial example of a stochastic fault-facies geocellular model, we consider a model based on the work presented in~\cite{Bjornaraa_etal_21}. 
The geometric distribution of facies is assumed to only vary in the vertical direction, while the facies flow properties can be expressed in terms of the shale gouge ratio (SGR), which is constant within each facies. The range of SGR values is confined to the interval $[0, 100]$, where the endpoints correspond to pure sandstone (SGR=0) and shale (SGR=100).
The model assumes that the total fault height $H$ and the number of facies $N_{\textup{facies}}$ is known.
The interior boundaries between the facies are assumed to be uniformly distributed around expected locations that are equidistant on $[0,\ H]$.
The assumption of no variation in the horizontal direction means that a realization of the fault geometry can be expressed in terms of a one-dimensional model.
The SGR distribution of a given facies is assumed to be independent of the other facies and to follow a Gaussian distribution with mean given by interpolation of the depth-dependent synthetic data shown in Figure~\ref{fig:SGR_data}, and an assumed standard deviation of 14 percentage points~\cite{Freeman_etal_08,Bjornaraa_etal_21}. The probability of this model yielding results outside the interval $[0, 100]$ is assumed to be negligible.
In total, the model has $2N_{\textup{facies}}-1$ random variables as inputs.

For each facies $j$ ($j=1,\dots, N_{\textup{facies}}$), there are various functional relationships to obtain permeability $k_j$ from the corresponding SGR, denoted $\text{SGR}_j$. For instance, in~\cite{Bjornaraa_etal_21}, where the permeability within facies $j$ is assumed constant, the relationship
\begin{equation}
\label{eq:SGR-perm}
\log(k_j) = 0.01\cdot\text{SGR}_j \cdot \log(k_{\text{c}}/k_{\text{s}}) + \log(k_{\text{s}}),  
\end{equation}
was given, where $\perm[\text{c}]$ is the permeability of clay-rich shale, and $\perm[\text{s}]$ is the permeability of sandstone. Other functional relationships between SGR and permeability are summarized in~\cite{Jolley_etal_07}, and a general discussion of methods for approximating the effect of clay content on fault permeability is provided in~\cite{Fisher_Jolley_07}.
In this work, we make the modeling decision to consider $\perm[\text{c}]=1$ mD as in~\cite{Bjornaraa_etal_21}, as well as $\perm[\text{c}]=0.001$ mD and $\perm[\text{c}]=0.0001$ mD based on the observation that many SGR-permeability models in the literature yield substantially lower output values than Eq.~\eqref{eq:SGR-perm} with $\perm[\text{c}]=1$ mD. 
For other SGR models that have been tested, the results aligned well with the stochastic modelling framework to be presented in this work, in the sense that the models themselves were different but of similar complexity and similar goodness-of-fit.

The upscaled uncertain permeability $K$ is computed from the weighted harmonic mean,
\begin{equation}
\label{eq:perm_upsc_rand}
K = \frac{\sum_{j=1}^{N_{\textup{facies}}} h_{j} }{\sum_{j=1}^{N_{\textup{facies}}}\frac{h_j}{\perm{j}}},
\end{equation}
where $h_j, k_j$ are the random height and permeability of facies $j$ for $j=1,...,N_{\textup{facies}}$. 
 In the case of single-phase flow, the permeability can be represented by a single random variable that does not depend on any other random variable.
To enable conditional sampling of the permeability (e.g., random sampling restricted to a subrange of the random variable), an empirical or parametric distribution (e.g., beta or lognormal) can be fitted to a large set of permeability samples. Conditional samples are then obtained from $F_{\text{fit.dist.}}^{-1}(u)$, where $F_{\text{fit.dist.}}^{-1}$ denotes the inverse cumulative distribution function of the fitted distribution, and $u$ is drawn uniformly from a subset of $[0,\ 1]$.

\begin{figure}[H]
\centering
\subfigure[Depth-dependent SGR.]
{\includegraphics[width=0.235\textwidth]{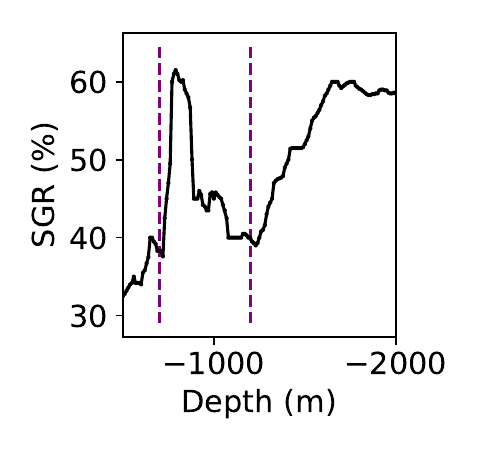}
\label{fig:SGR_data}
}
\subfigure[Permeability distribution, $k_{\text{c}}=10^{-4}$~mD.]
{\includegraphics[width=0.235\textwidth]{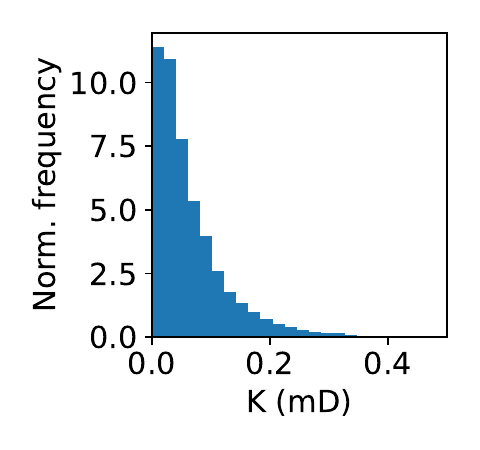}
\label{fig:perm_hist_k4}
}
\subfigure[Permeability distribution, $k_{\text{c}}=10^{-3}$~mD.]
{\includegraphics[width=0.235\textwidth]{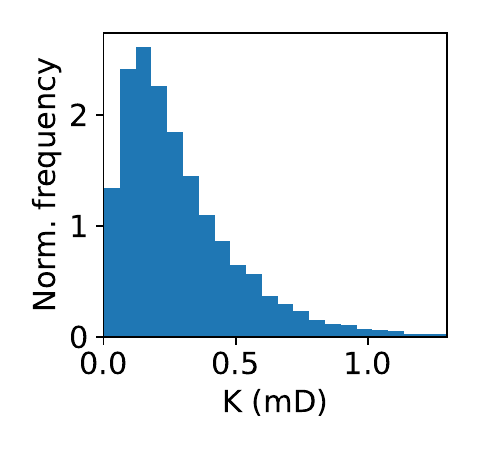}
\label{fig:perm_hist_k3}
}
\subfigure[Permeability distribution, $k_{\text{c}}=1$~mD.]
{\includegraphics[width=0.235\textwidth]{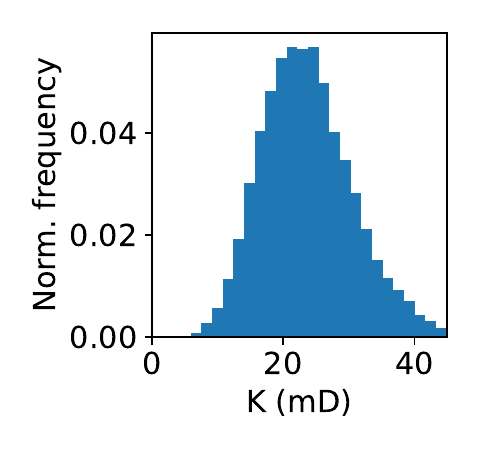}
\label{fig:perm_hist_k1}
}
\caption{Mean SGR based on synthetic data and upscaled permeability distributions from Eq.~\eqref{eq:SGR-perm}, \eqref{eq:perm_upsc_rand}, and different clay permeabilities. SGR assumed to follow a Gaussian distribution with a standard deviation 14\%.} 
\label{fig:SGR_dist_2}
\end{figure}

Empirical histograms, approximating the probability density functions (PDFs) for the permeability obtained from~\eqref{eq:perm_upsc_rand} and the three SGR-permeability models are shown in Figures~\ref{fig:perm_hist_k4},~\ref{fig:perm_hist_k3} and~\ref{fig:perm_hist_k1}.
An alternative to deterministic and fixed values of $\perm[\text{c}]$ and $\perm[\text{s}]$ is to account for the uncertainties in the composition by modeling them as random variables, but this approach will not be further investigated here.

\subsection{Upscaling of Two-phase Flow Functions}
\label{sec:twophase_upscaling}
Similar to the absolute permeability, the fault relative permeability and capillary pressure functions will be stochastic due to their dependency on the stochastic fine-scale facies.
However, since the modeling of the upscaled quantities must capture stochastic functions rather than just a number, as was the case for the absolute permeability, the task is more involved.

For the computations presented herein, the fine-scale capillary pressure function is computed from a Brooks-Corey function \cite{brooks1965hydraulic}. Assuming that the entry-pressure scales with the square root of permeability~\cite{leverett1941capillary} the entry-pressure in a geocell can be computed from the permeability in the geocell $k_{\text{cell}}$ and the entry-pressure $p_{ \text{s,entry}}$ and permeability $k_{\text{s}}$ in the sand as 
\[
p_{\text{cell}, \text{entry}} = p_{ \text{s,entry}} \cdot \sqrt{k_{\text{s}}/k_{\text{cell}}}.
\]
Thus, as for the absolute permeability, the fine-scale capillary pressure is a stochastic function of the SGR. The relative permeabilities are not stochastic on the fine scale, but will be so on the coarse scale due to the dependency of the fine-scale geometry and permeability distribution. For the fine-scale simulations herein, we use $p_{\text{s}, \text{entry}} = 2.5$ kPa and Brooks-Corey exponent of 0.67. These are values computed to fit the capillary-pressure values in the Smeaheia data-set (see section~\ref{sec:Smeaheia_model}). 

To compute upscaled two-phase flow functions for a given facies distribution, here we assume that the phases are in capillary equilibrium.
While this is a strong assumption that may be of questionable validity, e.g., in cases where viscous forces have significant impact on the flow, it allows for analytical upscaling of relative permeabilities and capillary pressure.
Specifically, we follow the approach described in~\cite{rabinovich2016analytical}, which for completeness is summarized in a form tailored to the setting considered here in~\ref{app:upscaling}.
The result is a set of discrete values for the flow functions,  in the sense of point values for permeabilities, and the corresponding values for capillary pressures, saturations, and relative permeabilities parametrized via the deterministic parameter $\sdisc$.
This is related, but not equivalent, to a fine-scale saturation at which the flow functions are sampled; see the appendix for an explanation.

\subsection{Reduced Multivariate Stochastic  Model for Flow Functions}
\label{sec:reduced_stochastic_model}
We apply the upscaling described in section~\ref{sec:twophase_upscaling} to a set of realizations of the fine-scale fault geometry, using the same realizations as those used to upscale permeability. Five coarse-scale flow functions are needed as input to the model $M$ for the case of two-phase flow: permeability $K$, capillary pressure $\upscaledPc$, saturation $\upscaledsat$, wetting relative permeability $K_{\text{r}}^{w}$, and non-wetting relative permeability $K_{\text{r}}^{\text{nw}}$.
 A large number of random samples of the  coarse-scale flow functions  evaluated at various deterministic fine-scale saturation values $\sdisc$ can be generated at negligible computational cost. 
The random flow functions are so far parameterized via the $2N_{\textup{facies}}-1$ random variables, but the upscaling implies that the stochastic dimensionality may be reduced to obtain a more succinct representation which is anyhow the format effectively accepted by the coarse-scale simulator. 

The stochastic representation of coarse-scale properties deviates from what would be expected of two-phase flow functions in a traditional deterministic modeling framework.
First, the model expresses the coarse saturation as a function of $\sdisc$ on the same level as the capillary pressure and relative permeabilities. This is contrary to deterministic modeling, where the saturation is considered an independent variable, and the capillary pressure and relative permeability are functions of this saturation.
This is simply a consequence of the combination of the fine-scale stochastic description and the upscaling procedure; as will be seen, the coarse-scale flow functions have the expected dependency of saturation, albeit in a stochastic sense, but deriving this relation is not straightforward.
Second, although the absolute permeability is not a deterministic function of saturation, nor of the two-phase flow functions, there is a clear random dependence between these parameters, as will be shown below.

\begin{figure}[H]
\centering
\subfigure[Upscaled data, $k_{\text{c}} =0.0001$ mD.]
{\includegraphics[width=0.99\textwidth]{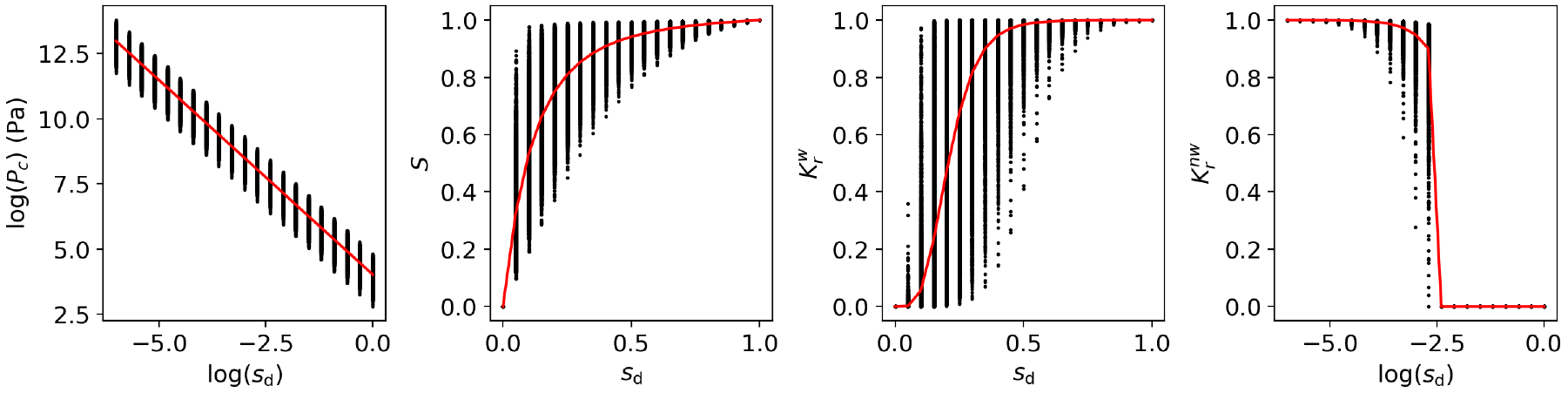}
}
\subfigure[Upscaled data, $k_{\text{c}}=0.001$ mD.]
{\includegraphics[width=0.99\textwidth]
{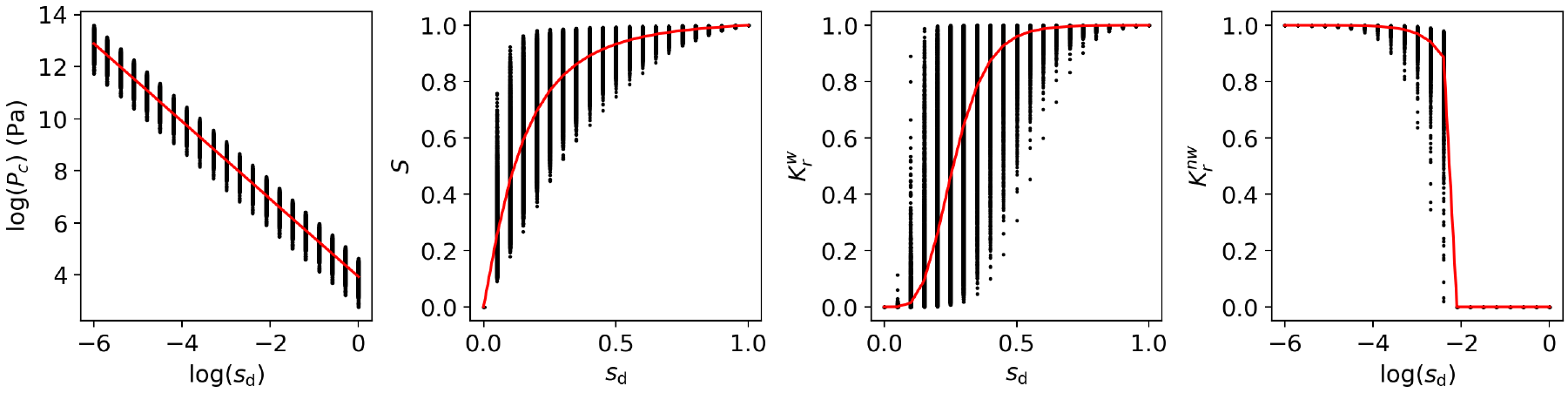}}
\subfigure[Upscaled data, $k_{\text{c}}=1$ mD.]
{\includegraphics[width=0.99\textwidth]
{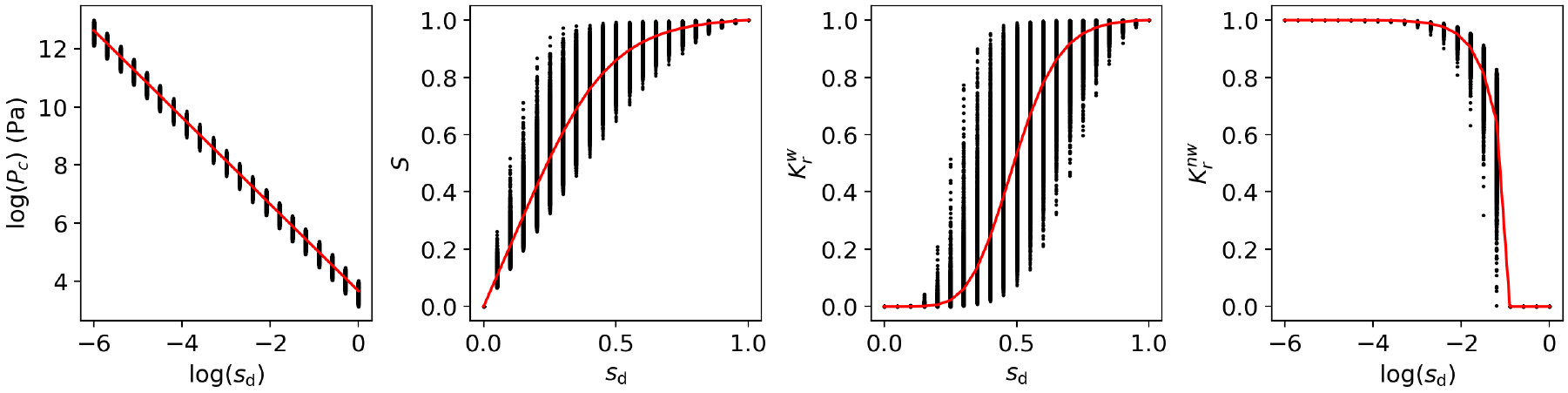}}
\caption{Visualization of random data (black dots) as support for stochastic modeling choices. Red curves indicate the means of the data. For ease of visualization, equidistant discretization of $\sdisc$ is used for capillary pressure and non-wetting relative permeability, and a log-scale discretization of $\sdisc$ for saturation and wetting permeability.}
\label{fig:param_param}
\end{figure}

Without any further conditions, the appropriate a priori model of the parameter quintuple $K$, $\upscaledPc, \upscaledsat, K_{\text{r}}^{w}, K_{\text{r}}^{\text{nw}}$ would be a five-dimensional stochastic process in $\sdisc$, for instance parameterized by a vector-valued Karhunen-Loeve expansion with suitable covariance structure~\cite{Perrin_13}. Such a representation would, however, still be complex and with a substantial computational cost and may not guarantee physically relevant output. These considerations therefore motivate us to search for a different method.
As we will demonstrate below, the data are indeed well represented by a multivariate random variable model with five distinct but dependent random variables, whose distributions are to be determined to match data at some user-defined value of $\sdisc$, henceforth denoted $\sdisc^{0}$. We then seek functional relationships to make sure the model matches the data for all relevant values of $\sdisc$. Ultimately, such a  model is fully and sufficiently motivated by presenting an acceptable goodness-of-fit and showing that it cannot be further simplified. Nevertheless, we find it instructive to illustrate the model choice via the data we seek to fit as follows. 
To that end, we generate $N_{\text{ref}}=10,000$ samples of the coarse-scale flow functions evaluated on each of 21  values of $\sdisc$ on the unit interval. Note that we have adapted the discretization to the scale of the $\sdisc$-axis, but in all numerical experiments we stick to the logaritmically spaced discretization that is denser for small values of $\sdisc$. The computation is done for three different values of \perm[c], and the results are shown in Figure~\ref{fig:param_param}.

The upscaled permeability is not a function of $\sdisc$, but so are the remaining four coarse-scale functions.
We first seek to characterize the dependence on the deterministic parameter $\sdisc$, and subsequently the residual is parameterized by random variables.
As can be seen from Fig.~\ref{fig:param_param}, for all values of \perm[c] there is a very good log-linear fit of the mean of the capillary pressure with respect to $\sdisc$. We note that the values taken by the capillary pressure at any two different values of $\sdisc$ are perfectly correlated, i.e,. the same single random parameter models all randomness in the capillary pressure.

The upscaled saturation and wetting relative permeability share some characteristics that make the stochastic modeling challenging. The randomness in each of them exhibits positive but not perfect two-point correlation (where perfect correlation means that the random outcome at any given value of $\sdisc$ completely determines the value at any other value of $\sdisc$). This lack of perfect correlation is due to crossing of sample paths, where a sample path here is to be understood as the flow function traced out by keeping the random variable it is dependent upon at a fixed value, and varying $\sdisc$ only. To be more specific, saturation sample path crossing means that there are two different values of $\sdisc$, say $\sdisc=s'$ and $\sdisc=s''$ (i.e., with $s'\neq s''$), and two distinct random samples $Y^{(m)} \neq Y^{(n)}$, for which we have for the saturation $\upscaledsat$ that $\upscaledsat(s', Y^{(n)}) > \upscaledsat(s', Y^{(m)})$ and $\upscaledsat(s'', Y^{(n)}) < \upscaledsat(s'', Y^{(m)})$. A suitable parametric stochastic model should reflect the correlation while ensuring that the output is monotone in $\sdisc$ to be compatible with the physical interpretation, and within physically relevant bounds, e.g., the unit interval for saturation and relative permeabilities.
Finally, the mean of the upscaled non-wetting relative permeability $K_{\text{r}}^{\text{nw}}$ can be approximated by a piecewise linear function in $\sdisc$. There is nearly perfect correlation in the randomness of this function, and it can be parameterized with a single random variable and a small correction for unphysical values that occur with small probability.
These properties would become more clear by plotting $K_{\text{r}}^{\text{nw}}$ as a function of $\sdisc$, but we have chosen to plot $K_{\text{r}}^{\text{nw}}$ as a function of $\log(\sdisc)$ in Figure~\ref{fig:param_param} for visibility of the full range of $\sdisc$.

Based on the above considerations and the observed data, an appropriate stochastic model for the upscaled parameters is given by

\begin{align}
& K = \exp(Y_1), \label{eq:K_Y1} \\
& \log(\upscaledPc(\sdisc)) = a_{\upscaledPc} + b_{\upscaledPc}\log(\sdisc) + Y_2, \label{eq:Pc_Y2}\\
& \log(\upscaledsat(\sdisc)) =  F^{-1}_{Y_{3,\sdisc}} \left( T^{\upscaledsat}_{\sdisc} F_{Y_3}(Y_3) \right), \label{eq:S_Y3}\\
& \log(K_{\text{r}}^{w}(\sdisc)) = 
F^{-1}_{Y_{4,\sdisc}} \left( T^{K_{\text{r}}^{w}}_{\sdisc} F_{Y_4}(Y_4) \right), \label{eq:Kw_Y4} \\
& K_{\text{r}}^{\text{nw}}(\sdisc) = 
\left\{
\begin{array}{ll}
 1+(a_{K_{\text{r}}^{\text{nw}}} + b_{K_{\text{r}}^{\text{nw}}} Y_5^{(\sdisc)}) \sdisc
 & \mbox{if } \sdisc \leq \sdisc^{0} \\
0 & \mbox{if } \sdisc > \sdisc^{0}
 \end{array}
 \right.
\label{eq:Knw_Y5}
\end{align}
with
\[
Y_5^{(\sdisc)} = \left\{
\begin{array}{ll}
Y_5 & \mbox{if } Y_5 \geq -(1+a_{K_{\text{r}}^{\text{nw}}}\sdisc)/(b_{K_{\text{r}}^{\text{nw}}}\sdisc) \\
-(1+a_{K_{\text{r}}^{\text{nw}}}\sdisc)/(b_{K_{\text{r}}^{\text{nw}}}\sdisc) & \mbox{if } Y_5 <- (1+a_{K_{\text{r}}^{\text{nw}}}\sdisc)/(b_{K_{\text{r}}^{\text{nw}}}\sdisc)
\end{array}
\right.
\]
where $\{ Y_j\}_{j=1}^{5}$ is a set of five dependent random variables, to be estimated from data. Due to the orders-of-magnitude variability of permeability an exponential relationship has been employed, although this is not necessary in a conceptual setting. In addition, we note that the log transformations in~\eqref{eq:S_Y3}-\eqref{eq:Kw_Y4} are not necessary, but they do make sure that the output values are positive.

Here, $F_{Y_3}$ (and $F_{Y_4}$) denotes the cumulative distribution function (CDF) of a given random variable $Y_3$ (and $Y_4$), where the dependence on $\sdisc = \sdisc^{0}$ has been suppressed for brevity of notation. The notation $T^{\upscaledsat}_{\sdisc}$ (and $T^{K_{\text{r}}^{\text{w}} }_{\sdisc}$) denotes a mapping from the unit interval to the unit interval, to account for correlation between data of $\upscaledsat$ (and $K_{\text{r}}^{\text{w}}$) at $\sdisc^{0}$ and other values of $\sdisc$. The purpose of the latter is to account for the fact that sample paths of $\upscaledsat$ or $K_{\text{r}}^{\text{w}}$ (as functions of $\sdisc$) can cross each other (i.e., a sample path of $\upscaledsat$ crosses another path of $\upscaledsat$). 
The mapping $T^{\upscaledsat}_{\sdisc}$ is constructed as follows, aiming to reflect the underlying stochastic dependencies. 
Recall that we have chosen a fixed $\sdisc = \sdisc^{0}$ as reference for the parameterization of the random variables $\{ Y_j \}_{j=1}^{5}$, i.e., we will use data for this particular value $\sdisc^{0}$ only to define and estimate the joint distributions of $\{ Y_j \}_{j=1}^{5}$. Any value of $\sdisc$ can be chosen for this purpose, but it is natural to choose the smallest considered $\sdisc$ as a starting point for a random process. The upscaled data at $\sdisc^{0}$ are then sorted to be used as a reference order. For any target $\sdisc > \sdisc^{0}$, the data are sorted locally with respect to the reference data at $\sdisc^{0}$. If no crossing of sample paths has occurred, the reference and target indices will be identically ordered, and $T^{\upscaledsat}_{\sdisc}$ would be the identity operator. In general, the order of indices will differ, and a one-to-one mapping between them is given by pairing entry by entry in the order of appearance. The range and domain of the mapping can both be normalized to the unit interval to reflect the fact that we want to map between the ranges of cumulative distribution functions. It is however practical to consider an integer-index discretization as follows. Let $\mathcal{I}$ be the indices permuting the reference data to the target data, and $\mathcal{I}^{-1}$ the indices that define the inverse mapping. We then introduce the approximation

\[
T^{\upscaledsat}_{\sdisc} F_{Y_3}(Y_3) \approx N_{\text{ref}}^{-1} \mathcal{I}^{-1}[N_{\text{ref}} F_{Y_3}(Y_3)],
\]
where double brackets denote integer rounding. Thus, the model $T^{\upscaledsat}_{\sdisc}$ approximates the two-point correlation present in the data,  which is not stationary, i.e., the correlation depends on the actual values of $\sdisc$, and not just the difference between any two points.

The mapping $F^{-1}_{Y_{3,\sdisc}}$ is given by the inverse of the empirical CDF at any given value of $\sdisc$, and ensures that the pointwise distribution approximates that of the underlying data. Hence, the accuracy of both the mappings $T^{\upscaledsat}_{\sdisc}$ and $F^{-1}_{Y_{3,\sdisc}}$ grows with the number of reference samples $N_{\text{ref}}$.
The discretization of the model for $K_{\text{r}}^{w}$ is completely analogous to that of $\upscaledsat$.

To summarize the above considerations, the stochastic facies description and upscaling method applied here, results in an upscaled stochastic model that can be described by relatively simple functional forms. 
We emphasize that changing either the facies description or the upscaling procedure may necessitate a new analysis of the upscaled model, and the stochastic model may take a different form, but the basic modeling paradigm and steps in the calculation remain the same.

\subsection{Copula Models and Transformations for Multivariate Distributions}
\label{sec:copula}
The stochastic upscaled model \eqref{eq:K_Y1}-\eqref{eq:Knw_Y5} provides a probabilistic description of the fault flow function in terms of coarse-scale parameters.
However, the random variables $\{ Y_{j} \}_{j=1}^{5}$ do not follow known distributions, nor are they independent. Both conditions are necessary for the model to be amenable to the computational uncertainty propagation method applied below.
We therefore parameterize and transform our stochastic model, using a combination of two techniques: First, copulas are employed to derive a stochastic model which honors the data, and have the attractive features of separation of dependency structure and marginal distributions, and a parameterization of (albeit dependent) uniform distributions on the unit interval.
Second, dependencies between the uniform variables are handled by transforming them to independent uniform variables using the inverse Rosenblatt transform~\cite{Rosenblatt_52}. Thus, the model provides an accurate mapping between standard and independent uniform random variables suitable for uncertainty propagation methods, and physically accurate random flow functions in the correct format for input to a reservoir simulator.
The application of both copulas and Rosenblatt transformations are described in some detail below, starting with an outline of copula construction of multivariate distributions.

According to Sklar's theorem~\cite{Sklar_59}, any multivariate CDF $F_{Y}(y_1,\dots, y_n)$ of the random vector $Y=(Y_1,\dots, Y_n)$ with marginal CDFs $F_{Y_1}(y_1),\dots, F_{Y_n}(y_n)$ can be expressed as
\begin{equation}
\label{eq:CDF_via_Sklar}
F_{Y}(y_1,\dots, y_n) = C_{Y}(F_{Y_1}(y_1),\dots, F_{Y_n}(y_n)),
\end{equation}
where $C_{Y}$ is a so-called $n$-copula. 
Differentiating the above expression and using the notation $c_{Y}(u_1,\dots, u_n) \equiv \partial_{u_1}\dots \partial_{u_n}C_{Y}(u_1,\dots,u_n)$, the joint PDF can be written,
\[
f_{Y}(y_1,\dots, y_n) = c_{Y}(F_{Y_1}(y_1),\dots, F_{Y_n}(y_n)) \prod_{j=1}^{n} f_{Y_j}(y_j),
\]
where $c_Y$ is the copula density function, and $f_{Y_j}(y_j)$ denote the marginal PDFs of the individual random variables $Y_1,\dots, Y_n$.
The copula framework separates the multivariate dependence properties described by the copula from the modeling of the marginal distributions. There are several commonly used families of copulas with different properties suitable for different kinds of data. Having chosen a suitable copula family, one can use data to estimate its parameters, and then use the fitted copula model to generate new samples from a multivariate distribution that approximates the distribution of the original data. 

Vine copulas offer a framework for multidimensional dependence modeling, relying on the use of repeated conditioning of up to two variables on subsets of all the remaining random variables. This leads to replacement of a multivariate function by a factorization into members of families of standard bivariate copulas~\cite{Bedford_Cooke_02}. The literature on high-dimensional copulas is more limited than that on bivariate copulas, and existing high-dimensional copulas can fail to correctly capture tail dependence, i.e., extreme behavior in several variables simultaneously. Conditioning on bivariate copulas also adds flexibility in the sense that different families of copulas can be used for different subsets of variables, allowing for a more accurate fit~\cite{Aas_Berg_09}.
An example of factorization by conditioning on bivariate copulas is provided in Appendix~\ref{sec:cop-fact-example}.
Note that the copula factorizations are exact, but in practice each bivariate distribution is estimated from data, introducing an approximation error.
Using data, both an appropriate vine copula factorization, bivariate copulas, and marginal distributions can be estimated. Pair-copula modeling, to which the class of Vine copulas belong, are suitable for high-dimensional modeling with the number of pair-copulas to be estimated being equal to $n(n+1)/2$~\cite{Aas_Berg_09}. However, this can become computationally prohibitive if the number of degrees of freedom from a large number of candidate models become very large~\cite{Brechmann_Schepsmeier_13}.

Once the copula model has been fitted to data, it can be sampled to provide inputs for uncertainty propagation through a physical model.
The random variables $\hat{U}_{j} = F_{Y_j}(Y_j)$ are uniformly distributed on the unit interval $(0,1)$, but they are, in general, not independent. For efficient coupling with many uncertainty propagation methods, a mapping between $\{ \hat{U}_{j} \}_{j=1}^{n}$ and a set of independent uniformly distributed random variables $\{ U_{j} \}_{j=1}^{n}$ is required. Such a mapping is offered by the inverse of the Rosenblatt transformation~\cite{Rosenblatt_52}, from the conditional variables $\{ U_{j}\}$ to the marginals $\{ \hat{U}_{j} \}$ which are the arguments of the copula~\eqref{eq:CDF_via_Sklar}. Specifically, we have the two sets of random variables:
\begin{align*} 
    \hat{U}_1 &= F_{Y_1}(Y_1), & \quad  U_1 &= F_{Y_1}(Y_1),\\
    \hat{U}_2 &= F_{Y_2}(Y_2), & \quad  U_2 &=F_{Y_2 | Y_1}(Y_2 | Y_1),\\
     & \vdots & \qquad &\vdots \\
    \hat{U}_n &= F_{Y_n}(Y_n), & \quad  U_n &= F_{Y_n | Y_1,\dots, Y_{n-1}}(Y_n|Y_1,\dots, Y_{n-1}).
\end{align*}
The output from the inverse Rosenblatt transformation then gives a means to express the target random variables $\{ Y_j \}_{j=1}^{n}$ as a function of the parametric input distributions:
\begin{align*} 
    Y_1 &= F^{-1}_{Y_1}(U_1)\\
    Y_2 &= F^{-1}_{Y_2 | Y_1}(U_2 | Y_1)\\
     & \vdots  \\
    Y_n   &= F^{-1}_{Y_n | Y_1,\dots, Y_{n-1}}(U_n|Y_1,\dots, Y_{n-1}),
\end{align*}

\subsection{Results: Copula Fit  for Upscaled Flow Functions}
\label{sec:copula_fit}
The upscaling model described in previous sections is used to generate $N_{\text{ref}}=10,000$ samples of the dependent random variables $\{ Y_j\}_{j=1}^{5}$ given by Eqs.~\eqref{eq:K_Y1}--\eqref{eq:Knw_Y5}, and also as shown in Figure~\ref{fig:param_param}. These samples are treated as data to which a copula model is fitted to allow targeted sampling and application of various uncertainty propagation methods. The Python interface pyvinecopulib to the C++ vinecopulib package~\cite{Nagler_etal_22package, Czado_Nagler_22} is employed for this objective, resulting in copula models defined by the tree structures illustrated in Figure~\ref{fig:vine-trees} in Appendix~\ref{sec:cop-tree-fit}.

\begin{figure}[H]
\centering
\subfigure[$k_{c}=0.0001$ mD.]
{\includegraphics[width=0.99\textwidth]{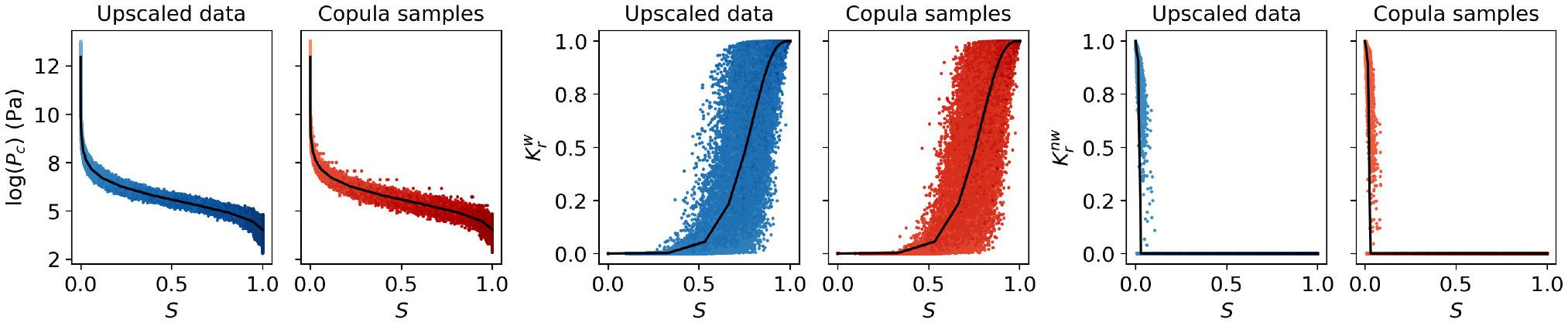}}
\subfigure[$k_{c}=0.001$ mD.]
{\includegraphics[width=0.99\textwidth]{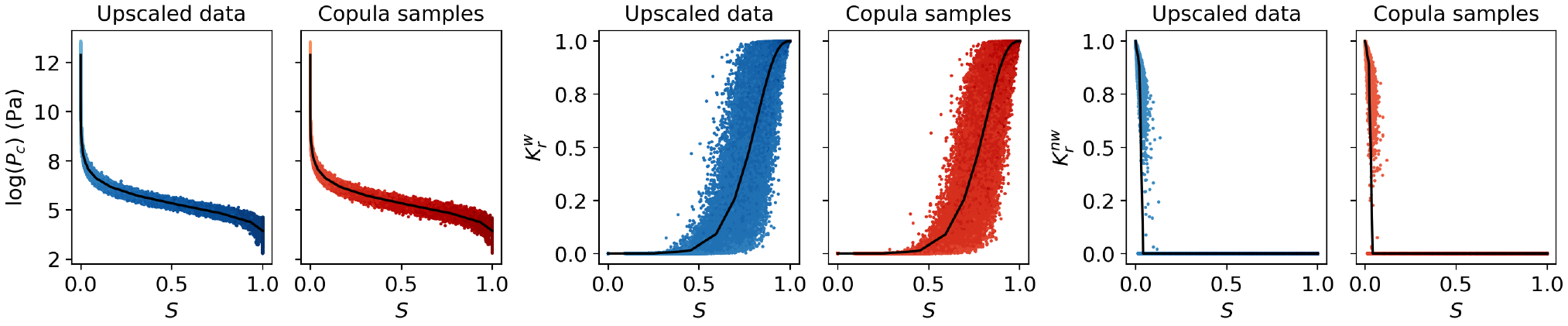}}
\subfigure[$k_{c}=1$ mD.]
{\includegraphics[width=0.99\textwidth]{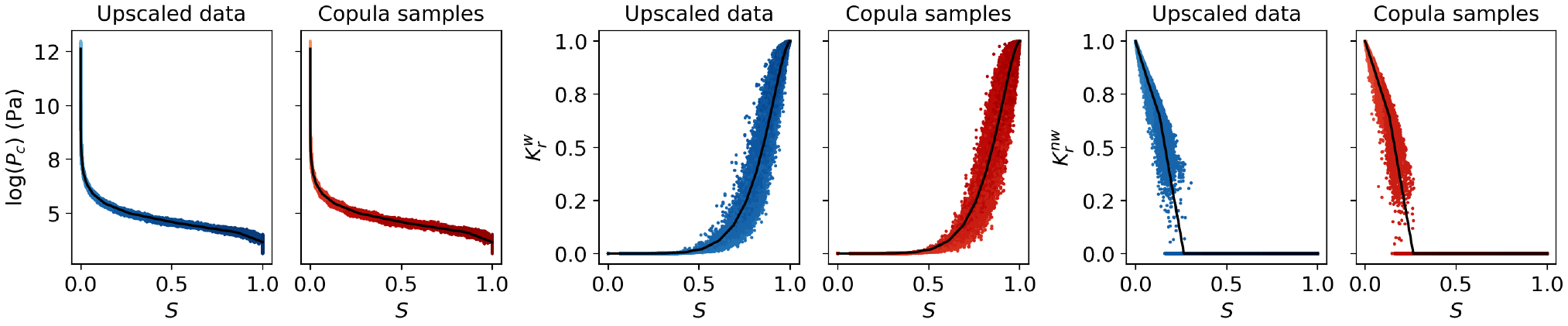}
}
\caption{Upscaled and copula models of capillary pressure  and relative permeabilities as function of saturation. The coloring of the data runs from low deterministic entry pressure (light blue/red) to high entry pressure (dark blue/red). The mean values of the data are shown in black.} 
\label{fig:Bjornaraa_model_both_perms_subset}
\end{figure}

It is of interest to illustrate the effect of representing the upscaled model in the form of copulas, and also to show the final stochastic flow functions in a form that is comparable with traditional deterministic modeling.
To that end, Figure~\ref{fig:Bjornaraa_model_both_perms_subset} plots the capillary pressure and relative permeability as functions of the saturation.
The figure shows both data upscaled from fine-scale realizations and data sampled from the copula model, using different values of clay permeability $\perm[\text{c}]$ in~\eqref{eq:SGR-perm}.
We emphasize that since the samples are drawn independently from the two distributions (upscaled data vs copula model), they will not exhibit a one-to-one relation. However, for a sufficiently large sample size, the generated point clouds should have a similar shape.
We observe that, although there are some discrepancies, there is overall a good agreement between the results from the fine-scale and the copula model.
Moreover, the plots indicate that the stochastic flow functions take forms that are reminiscent of what would be expected from a deterministic model.

\begin{figure}[H]
\centering
\subfigure[$k_{c}$=0.0001 mD.]
{\includegraphics[width=0.48\textwidth]{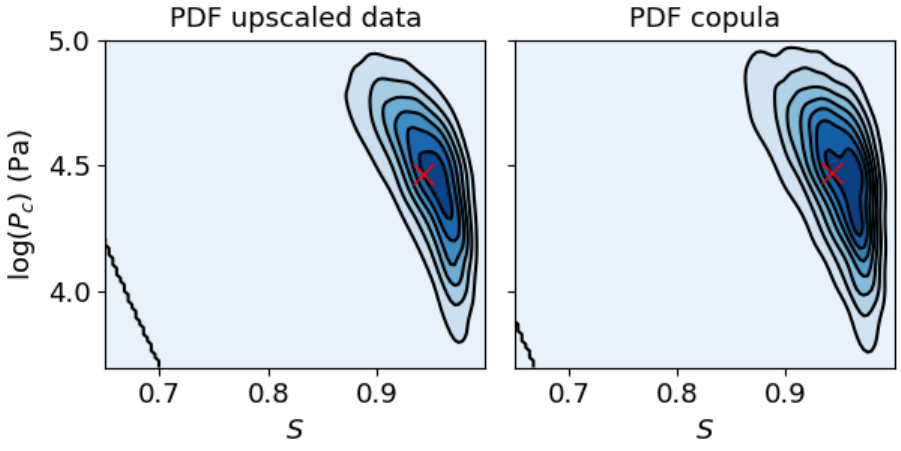}
}
\subfigure[$k_{c}$=0.0001 mD.]
{\includegraphics[width=0.48\textwidth]{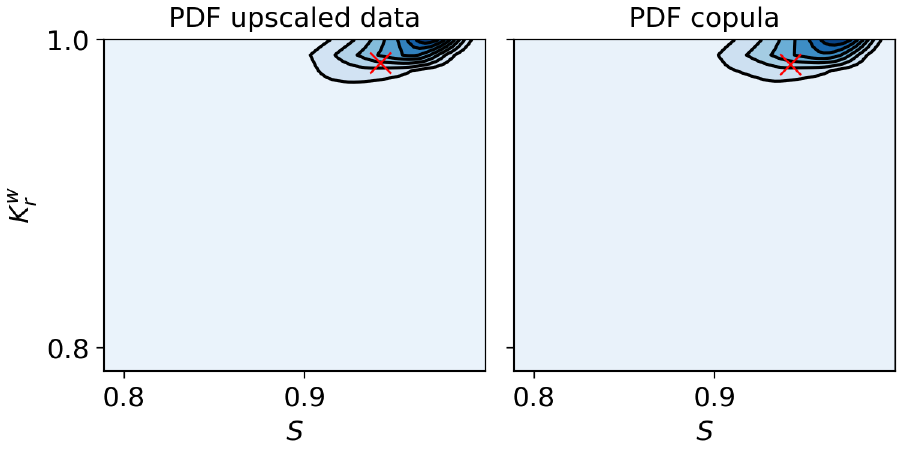}
}
\subfigure[$k_{c}$=0.001 mD.]
{\includegraphics[width=0.48\textwidth]{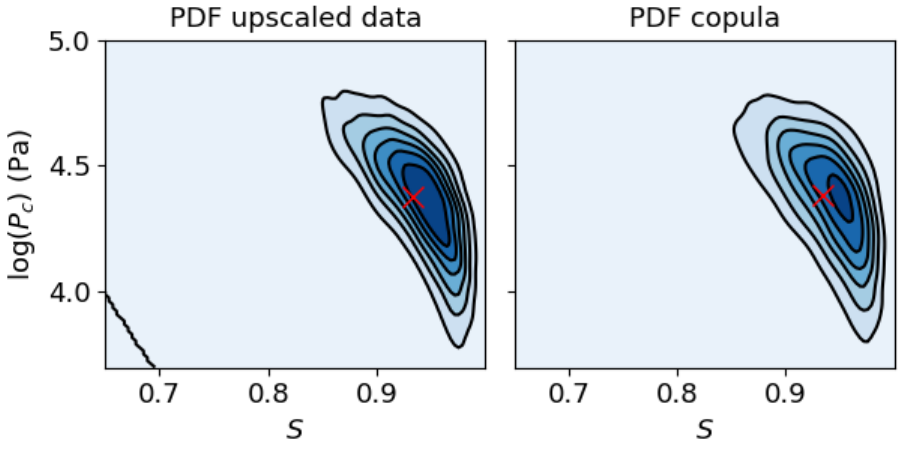}
}
\subfigure[$k_{c}$=0.001 mD.]
{\includegraphics[width=0.48\textwidth]{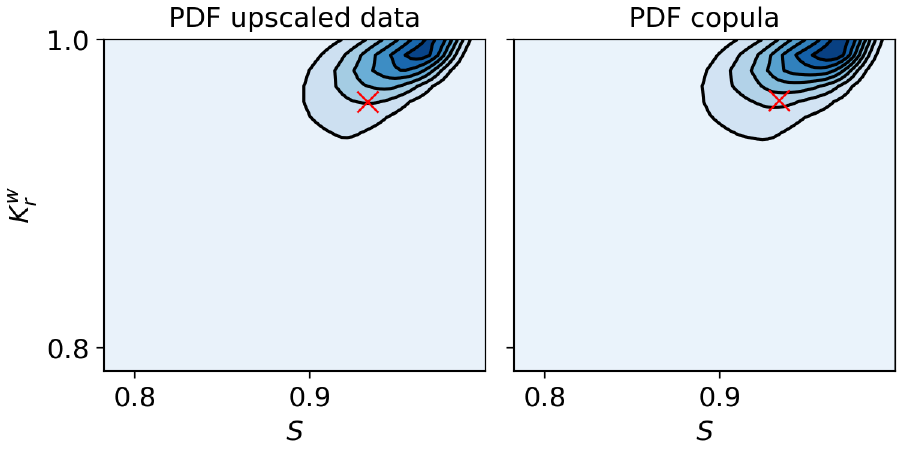}
}
\subfigure[$k_{c}$=1 mD.]
{\includegraphics[width=0.48\textwidth]{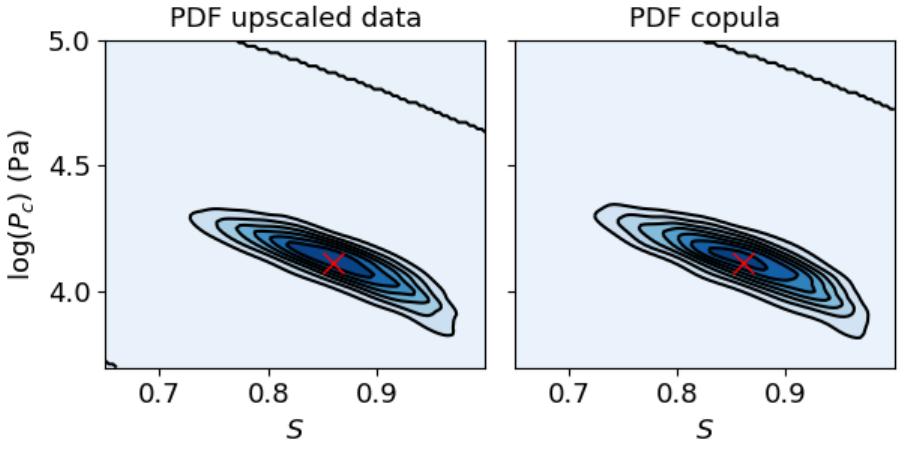}
}
\subfigure[$k_{c}$=1 mD.]
{\includegraphics[width=0.48\textwidth]{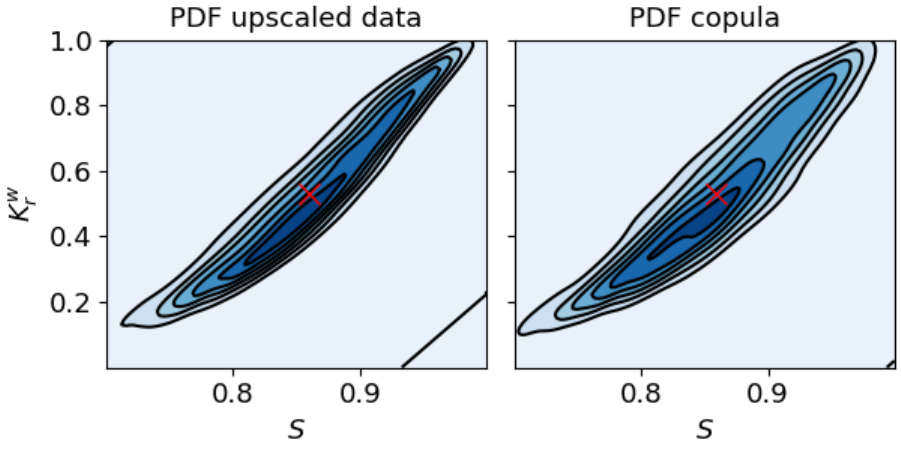}
}
\caption{
Estimated PDFs for $\sdisc=0.5$ for the joint marginal distribution of upscaled capillary pressure and saturation, and upscaled wetting relative permeability and saturation, respectively.
The PDFs are estimated using kernel density estimation from 10,000 samples. Expected values are marked with red crosses.} 
\label{fig:Bjornaraa_model_low_perm_local}
\end{figure}

As a second illustration and as a measure of similarity between the underlying and the modeled  distributions, we show estimated PDFs for the capillary pressure and wetting relative permeability, both as functions of saturation, for the three values of $\perm[c]$.  Contour plots of kernel density estimates of the PDFs of the upscaled data and equally many samples from the estimated copula and inverse Rosenblatt transformation of independent uniform samples are shown in Figure~\ref{fig:Bjornaraa_model_low_perm_local}. The figure axes are scaled to include the whole range of observed data, which can be sparse in certain regions, leading to a localization of the estimated PDFs to a small subregion of the figure.
Since the PDF is a point value, we consider a fixed value of $\sdisc=0.5$, and note that a different value would imply a different PDF. For this value, the non-wetting relative permeability is zero, and PDFs are shown only for capillary pressure and wetting relative permeability as functions of saturation.
Again, we observe a good agreement between the data generated from the fine-scale and the copula model.
Moreover, a comparison of the results for the three $\perm[\text{c}]$ values shows that, both the capillary pressure and wetting relative permeability vs saturation functions have PDFs with relatively similar shape for the two smaller values of $\perm[\text{c}]$. For $\perm[\text{c}]=1$ mD there are significant differences in the shape of the PDFs and in their expected values (marked with red crosses in Figure~\ref{fig:Bjornaraa_model_low_perm_local}).
In particular, the wetting relative permeability vs saturation exhibits marked differences in the shape of the PDF for the different fine-scale facies properties, in particular due to the proximity to the upper bound $K_{\text{r}}^{w}\leq 1$ for $\perm[\text{c}] = 0.001, 0.0001$ mD.

\begin{figure}
    \centering  \subfigure
{\includegraphics[width=0.99\textwidth]{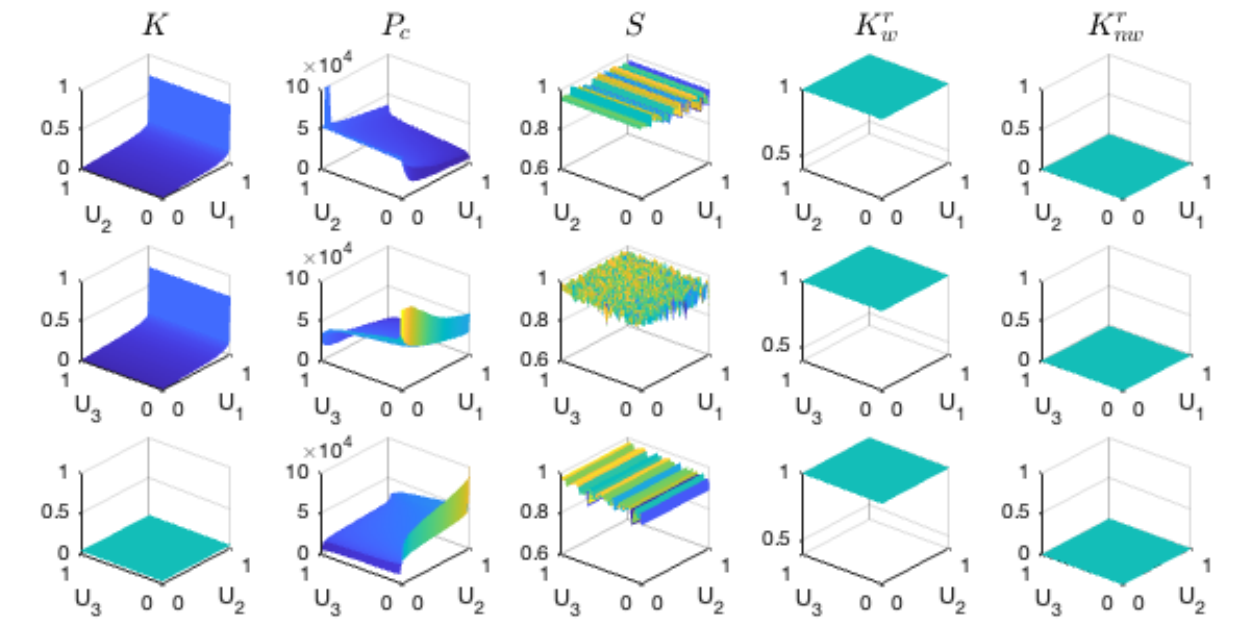}}
    \caption{Upscaled flow functions displayed as functions of pairs of $U_1,\dots,U_3$ for $\sdisc=0.5$, $k_{\text{c}}=0.0001$ mD.} 
    \label{fig:flow_funs_vs_Us}
\end{figure}

Figures~\ref{fig:param_param}, \ref{fig:Bjornaraa_model_both_perms_subset} and~\ref{fig:Bjornaraa_model_low_perm_local} indicate good fits with respect to data of the copula model. Many UP methods take independent (and often uniform) random variables as inputs, and the properties of the output as a function of those inputs will typically determine the performance of the selected UP method. Before presenting the UP method, it is therefore instructive to visualize the five upscaled model outputs as functions of $U_j$, $j=1,\dots, 5$. Figure~\ref{fig:flow_funs_vs_Us} depicts all upscaled flow functions of all combinations of pairs of the first three independent uniforms $\{U_j \}_{j=1}^{3}$, where the remaining random variables are kept fixed at 0.5. The remaining seven combinations of all pairs of $\{U_j \}_{j=1}^{5}$ are exhibiting similar patterns, and are therefore not included. The results are shown for $\perm[c]=0.0001$ mD and $\sdisc=0.5$ only. Other values of $\sdisc$ would give different plots, but some general observations will be useful. The upscaled non-wetting relative permeability is identically zero for this value of $\sdisc$ and hence displays no variance at all, whereas upscaled capillary pressure varies in all input variables. The wiggles in upscaled saturation are not spurious oscillations but a consequence of the nonlinear mapping $T^{\upscaledsat}_{\sdisc}$ in Eq.~\eqref{eq:S_Y3}. Although these qualitative observations do not make it apriori clear what an optimal UP method would be, they do indicate that a method that is insensitive to irregular features and that can exploit localization in random space may be a suitable choice for the problems of interest.  Figure~\ref{fig:flow_funs_vs_Us} indicates that there will be large regions of the input parameter space where the flow functions do not vary, and hence there would be no variation in an output quantity of interest. This can be exploited by adaptive sampling, placing fewer points in low-variability regions in favor of denser sampling in regions with large variability. Adaptive stratified sampling, to be described next, does not rely on regularity assumptions and is tailored to problems where adaptivity in stochastic space is possible, i.e., with localization of variance. Hence, we believe it is an appropriate choice of method for the problems considered in the current work.

\section{Adaptive Stratified Sampling for Uncertainty Propagation}
\label{sec:adaptive_sampling}
We now briefly outline adaptive stratified sampling and refer the reader to~\cite{Pettersson_Krumscheid_22} for more details. Recalling Eq.~\eqref{eq:high_level_form} and suppressing the dependence on physical deterministic parameters $\eta$ for ease of notation, let $Q=M(Y(U)) \equiv \tilde{M}(U)$ denote the quantity of interest as a function of random input variables $U \in [0,\ 1]^{n}$ where we have assumed that all inputs have been transformed to the unit hypercube, for instance by means of the Rosenblatt transformation.  In standard Monte Carlo (SMC) sampling, the estimator of the expectation of a quantity of interest is approximated as
\begin{equation}
\label{eq:smc_est}
\mu_{Q} \equiv \E{Q}=\E{\tilde{M}(U)} \approx \hat{Q}_{\text{MC}} \equiv \frac{1}{N}\sum_{j=1}^{N} \tilde{M}(u^{(j)}),
\end{equation}
where $\{ u^{(j)} \}_{j=1}^{N}$ denotes a sample set of size $N$ of standard pseudo-random samples in $[0,\ 1]^{n}$.
The variance of the estimator~\eqref{eq:smc_est} can be reduced by partitioning the stochastic domain into disjoint strata, and sampling each stratum independently according to a suitable allocation rule to determine how the sample budget will be distributed between strata. Common allocation rules include proportional allocation, where the number of samples per stratum is proportional to the size of the stratum, and optimal allocation, that minimizes variance if the intra-stratum variances are known. A hybrid between proportional and optimal allocation provides a good compromise between robustness and optimality, and will be employed here.

Let $\mathcal{S}$ denote the stratification (partition) of the stochastic domain $[0,\ 1]^{n}$, with disjoint strata of probability measure $p_{S}$ for all $S\in \mathcal{S}$. The probability measure or size is allowed to vary between strata. Then, if $N_S^{\text{prop}}$ and $N_S^{\text{opt}}$ denote respectively the number of samples allocated to stratum $S$ with proportional and optimal sampling, hybrid allocation corresponds to the allocation rule 
\begin{equation}
    N_S^{\alpha} := (1-\alpha) N_S^{\text{prop}} + \alpha N_S^{\text{opt}}, 
  \label{eq:allo_rules:hybrid}
\end{equation}
for some choice of $\alpha \in [0,\ 1]$. The hybrid stratified sampling estimator is given by
  \begin{equation}
  \label{eq:ss_estimator:hybrid}
  \hat{Q}_{\alpha} := \sum_{S\in\mathcal{S}}  \frac{p_S}{N_{S}^{\alpha}}\sum_{j=1}^{N_{S}^{\alpha}}Q_{S}^{(j)},
\end{equation}
for all $1\le j\le N_S^{\alpha}$ and $S\in\mathcal{S}$.

The relative performance of two UP methods is measured by comparing their variances, assuming both methods are unbiased.
The variance of the estimator $ \hat{Q}_{\alpha}$ is given by
\begin{equation} \label{eq:ss_estimator:var:hybrid}
  \var{\hat{Q}_\alpha} =  \frac{1}{N}\sum_{S\in\mathcal{S}}\frac{p_S \sigma_S^2}{1 + \alpha (\overline{\sigma}_S-1)}\;,
\end{equation}
where $\overline{\sigma}_S \equiv \sigma_{S}/ \sum_{T \in \mathcal{S}} p_{T}\sigma_{T}$.
As a quantitative measure of the performance of adaptive stratified sampling relative to standard Monte Carlo, we define
\begin{equation}
\label{eq:speedup}
\text{Speedup} \equiv \frac{\var{\hat{Q}_{\text{MC}}} }{\var{\hat{Q}_\alpha} }.
\end{equation}

The variance~\eqref{eq:ss_estimator:var:hybrid} is only exact if the prescribed sampling rates~\eqref{eq:allo_rules:hybrid} are maintained, which may be difficult in practice.
If that is not the case, for instance as a consequence of too few samples to compensate for previous non-ideal allocation, the resulting estimates may be inaccurate. As a remedy, the denominator in~\eqref{eq:speedup} can be computed by simply computing the sample variance over repeated independent evaluations of the  stratified sampling estimator. This procedure is robust but also expensive and is only motivated for demonstration purposes and cases where the model of interest is relatively cheap to sample from. In this work we use the Python implementation of the adaptive stratified sampling method available at~\cite{Krumscheid_etal_package_23}.

\section{Physical Model}
\label{sec:simulator}
Here we give the governing equations and the numerical solution approach for the physical model $M$.
We further describe the data set used in the numerical experiments presented in the next section.

\subsection{Numerical Reservoir Simulation Model}
A simulator is used to compute the quantity of interest from a set of physical parameters and stochastic variables. The simulator solves a set of discretized equations for mass conservation for CO$_2$ 
and brine 
in a two-phase system of gas(g) and liquid(l) using a blackoil formulation
\begin{equation*}
\begin{aligned}
 \frac{\partial \phi (b_l S_l)}{\partial t}+\nabla\cdot(b_l\pmb{v}_b)&=0, \\
 \frac{\partial \phi (b_g S_g + r_s b_l S_l)}{\partial t}+\nabla\cdot(b_g\pmb{v}_g + r_s b_l \pmb{v}_l)&=0,
 \end{aligned}
\end{equation*}
where $\phi$ is the rock porosity, $r_s$ the dissolved gas-liquid ratio, and $b_{g, l}$ the inverse formation-volume factor (measuring the ratio between the bulk volumes of the phases occupied at surface and reservoir conditions). The saturation sum is $S_g+S_l=1$ and the volumetric flux for phase $\alpha = {g,l}$ is given by the Darcy equation; 
\begin{equation*}
\pmb{v}_{\alpha}=-K \frac{K_{r}^{\alpha}}{\mu_{\alpha}}(\nabla p_{\alpha}-\rho_{\alpha}\pmb{g}),
\end{equation*}
where $\rho_\alpha$ is the fluid density, $p_{\alpha}$ is phase pressure, $\pmb{g}$ the gravity, $\mu_\alpha$ the fluid viscosity. The phase pressures are related through the capillary pressure $\upscaledPc(S_l)=p_g-p_l$.

These equations are discretized on a grid using a two-point-flux approximation with upwinding for the spatial discretization, and backward Euler for the time discretization. We refer to \cite{Aziz_Settari_79, rasmussen2021open} for details. The flux between two grid cells ($a$ and $b$) can then be formulated in terms of a transmissibility $T_{a-b}$;
\begin{equation}
    u_{\alpha, a-b} = -T_{a-b} \lambda_{\alpha}^{\text{up}} (p_{\alpha,a} - p_{\alpha,b} - \rho_{\alpha}^{\text{avg}}g).
\end{equation}
Here the mobility $\lambda^{\text{up}}_{\alpha} = k_{r\alpha} / \mu_{\alpha}$ is evaluated in the upstream cell while the pressure is evaluated in the cell centers and an average value of the cells is used for the density $\rho_{\alpha}$. The transmissibility between cell $a$ and $b$ is given by the harmonic mean of the one-sided transmissibilities such that
\begin{equation}\label{eq:trans}
T_{a-b} = \frac{T_{a} T_{b}}{T_{a} + T_{b}},  \quad T_{a,b} = K_{a,b} A / D_{a-b},
\end{equation}
where $K_{a,b}$ is the permeability in cell $a$ and $b$, respectively, and $D_{a - b}$ is the distance from cell $a$ to the center of the interface between the cells and from cell $b$ to the center of the interface, respectively. 
The discretized equations are solved to an appropriate tolerance using a Newton-Raphson type method with iterative linear solvers. For further details on solvers and discretization, we refer to~\cite{rasmussen2021open} and the OPM manual~\cite{opmFlowManual}. 

In addition to the equations presented herein, we need a set of boundary conditions and sink/source terms. For reservoir simulators, the sink/source terms are typically wells. For details on how wells are modeled in OPM Flow we refer to the OPM manual \cite{opmFlowManual}. The default setup assumes no-flow boundaries, but flow through the boundaries can be included using various aquifer models, both based on simplified numerical models and analytical approximations. In our setup, we will use numerical aquifers to model the flow to aquifers connected through faults.  

In the simulation model, dynamics within the numerical aquifers are simplified, and the aquifers are represented with only a few grid cells in the model that can be connected to any number of cells in the reservoir model using non-neighbor connections (NNC). Note that the framework itself allows for adding wells and more cells to the numerical aquifers if needed. Properties of the numerical aquifers like porosity, volume, permeability and connectivity can be set independently of the physical grid, which gives significant flexibility in usage. A very large porosity value can for instance be assigned to increase the pore volume of the aquifers, effectively disallowing pressure build-up in the aquifer.

The transmissibility used for the NNCs from the reservoir to the aquifer is given by the harmonic mean of the half transmissibility computed from the reservoir side, $T_{\text{r}}$ and the transmissibility computed from the aquifer side, $T_{\text{aq}}$, as
\begin{equation}\label{eq:num_aqf_trans}
T_{\text{r}-\text{aq}} = \frac{T_{\text{r}} T_{\text{aq}}}{T_{\text{r}} + T_{\text{aq}}},  \quad T_{\text{aq}} = 2 K_{f} A_{f} / L_{f}. 
\end{equation}
The fault area, $A_{\text{f}}$, represents the area available for flow between the reservoir cell and the aquifer, and the length, $L_{\text{f}}$, represent the length between the reservoir and the aquifer, while $K_{\text{f}}$ is the effective permeability of the fault. 
Note that we here assume that the transmissibility computed for the aquifer side is given by the properties of the fault, which is valid as long as 
$K_{\text{f}} A_{\text{f}} / L_{\text{f}} \ll K_{\text{aq}} A_{\text{aq}} / L_{\text{aq}}$. For along-fault flow, this assumption is true for most range of permeability as $A_{\text{f}} \ll A_{\text{aq}}$. For cross-fault flow $A_{\text{f}} = A_{\text{aq}}$ and the assumption depends on $K_{\text{f}} \ll K_{\text{aq}}$. 

Since the fault is modeled implicitly by the aquifer connections, we use the saturation function computed for the fault in the aquifer setup. As the relative permeability is evaluated in the upwind direction~\cite{Aziz_Settari_79}, a second cell is added to the aquifer to assure the effect of the upscaled fault relative permeability is modeled in the setup. The aquifer thus consists of two cells, the first with pore volumes representing the fault core, and the second with a pore volume representing the aquifer itself. 

\subsection{Description of the Smeaheia Data-set}
\label{sec:Smeaheia_model}

The starting point for the physical model used in the simulations is the Smeaheia data-set available at CO2Share (https://co2datashare.org/). To reduce the computational cost, we define a sector model representing a part of the reservoir by removing grid cells far from the injection point and the fault.
This reduces the number of cells from 250,000 to 44,000, and leads to a simulation time of less than a minute for a single model realization of $M$. A typical grid cell in the model is 400 m times 400 m in the $x$ and $y$ direction, respectively, and between 1 m and 20 m in the $z$ direction. While this is very coarse, it is also typical for field-scale simulations. Using a coarse grid near the fault may introduce modeling error due to inaccurate representation of pressure gradients and upconing effects near the fault. Analytically derived model corrections or local grid refinement could remedy the modeling error \cite{kang2014analytical}, but will not be considered in our approach. The reservoir consists of six horizontal layers with constant permeabilities. Three highly permeable sand layers with horizontal permeability 1000 mD, 1000 mD and 850 mD, respectively that are separated by two less permeable shale layers with horizontal permeability 50 mD and a less permeable bottom layer with horizontal permeability 25 mD. The vertical permeability is 1/10 of the horizontal permeability.    

\begin{figure}[H]
    \centering
    \subfigure[Smeaheia section model with injector  (red arrow). Orange shows the part of the Vette fault connected to the top aquifer, while red shows  connection to Troll.]{\includegraphics[width=0.48\textwidth]{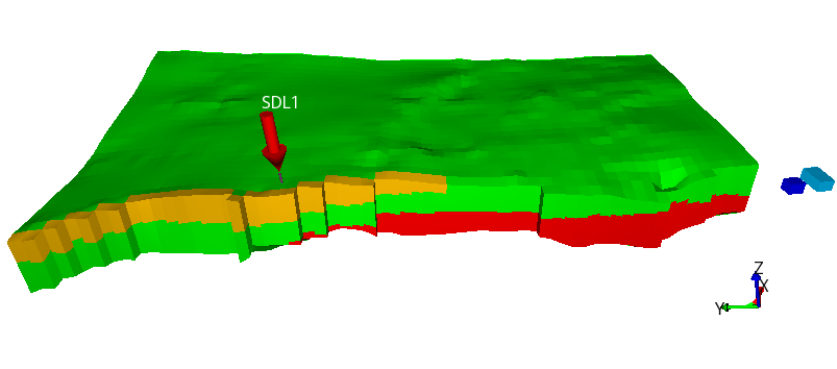}}
    \hspace{0.02\textwidth}
    \subfigure[Hydraulically connected zones are represented as numerical aquifers and coupled to the storage simulation via non-neighbor connections (NNC).]{\includegraphics[width=0.48\textwidth]{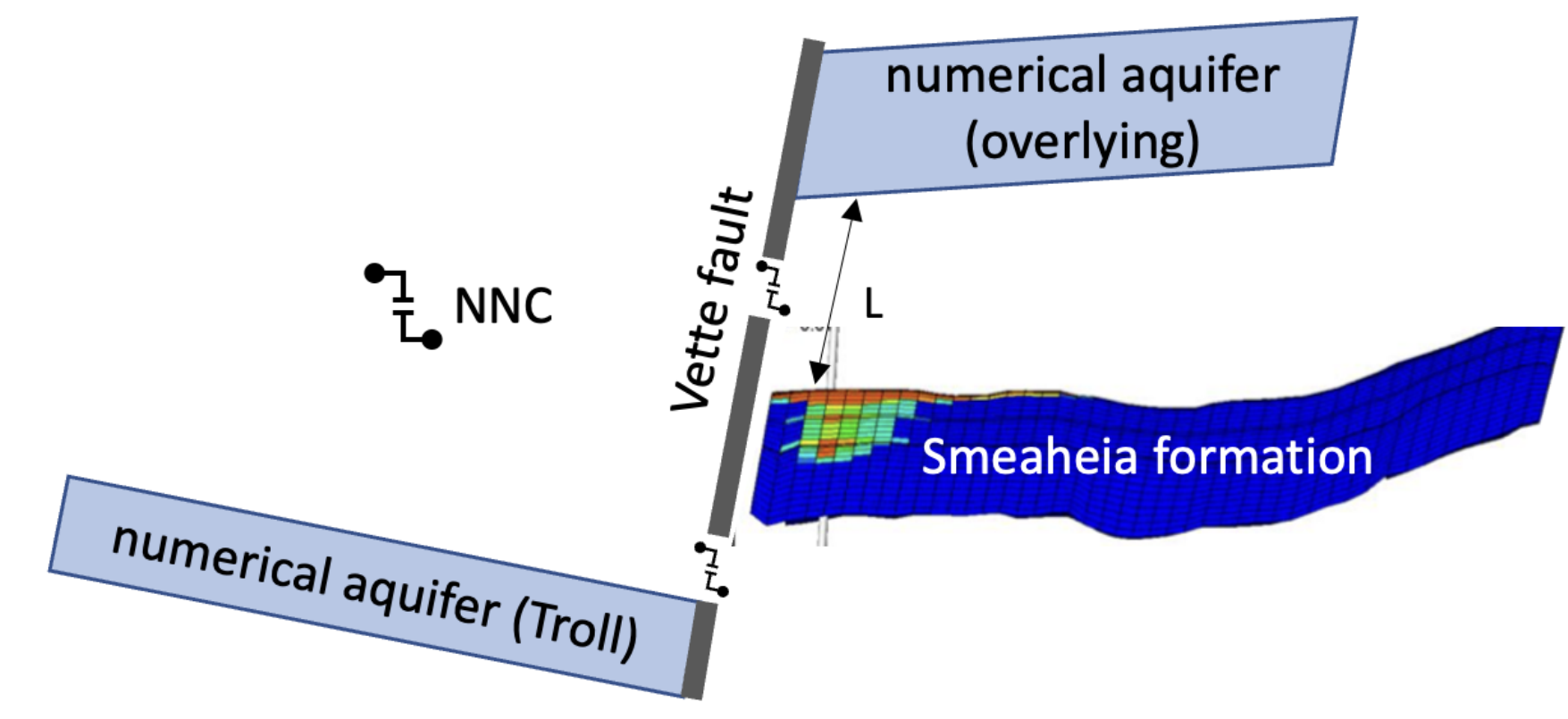}}
    \caption{The section of the Smeaheia simulation model used in the study and the setup with non-neighbor connections to hydraulically connected zones.
  }
    \label{fig:smeaheia}
\end{figure}

To include connections to nearby aquifers, the reservoir model is extended with two numerical aquifers, as illustrated in Figure~\ref{fig:smeaheia}. The first one, denoted Top aquifer, (shown in orange in Figure~\ref{fig:smeaheia}) represents an aquifer above the Smeaheia model potentially connected vertically along the Vette fault, while the second (shown in red in Figure~\ref{fig:smeaheia}) represents the Troll field connected through Vette in the lower parts of the fault and is denoted Troll aquifer. To represent the highly complex dynamics of the Troll field with only one cell is clearly a crude approximation, which nevertheless should suffice to investigate the impact of Troll on the Smeaheia formation through uncertain fault parameters.
An extra cell is added for the Top aquifer representing the fault in order to model the effect of having a different saturation function in the fault. For the Troll aquifer this is not needed as the connection is in the water zone and no two-phase effects are present.



There are two sets of saturation functions, i.e., relative permeability and capillary pressure, one for the reservoir and one for the fault. Each are discretized and provided to the model as input tables. The reservoir saturation functions are computed using a standard Brooks-Corey function. For the fault, the tabulated values are obtained from realizations of the stochastic upscaled functions as computed in sections~\ref{sec:twophase_upscaling}--\ref{sec:reduced_stochastic_model}, sampled with the adaptive stratified sampling approach described in section \ref{sec:adaptive_sampling}. 

The parameter set of the numerical aquifers are given in Table~\ref{tab:aquifer_input} and for details we refer to AQUNUM in~\cite{opmFlowManual}. CO$_2$ is injected with a constant rate of 1.6 million tons per year 
for 59 years. The injection rate is the rate given in the Smeaheia data-set. 
For simplicity, we do not consider CO$_2$ to dissolve into the brine phase. For the injection rate and timescale  
this simplification is not expected to alter the results qualitatively.

\begin{table}
    \centering
    \begin{adjustbox}{max width=\textwidth}
   \begin{tabular}{c|c c c c c c}
     Name & Area [m$^2$] ($A_{\text{f}}$) & Length [m] ($2L_{\text{f}}$) & Depth [m] & Poro [-] & Pres [bars] \\
     \hline
     Top Aquifer & 40 & $2\cdot250$ & * & 1e12 & * \\
     Troll Aquifer & 40 & $2\cdot250$ & 1100 & 1e12 & 110 
    \end{tabular}
    \end{adjustbox}
    \caption{Input numerical aquifers. The $*$ indicates that the values are defaulted and thus computed internally in the simulator.}
    \label{tab:aquifer_input}
\end{table}

\section{Numerical Results}
\label{sec:results}
Numerical results will be presented using the stochastic model described in section~\ref{sec:stochastic_faults} with the following setup.
A geocellular fine-scale model comprised of 20 geocells is employed for the Vette fault from 700 m to 1200 m depth. This model is represented by a total of 39 independent random parameters; 19 uniformly distributed uncertain geocell boundary locations, and 20 uncertain SGR values. Using the presented physical-stochastic upscaling framework results in coarse-scale two-phase flow parameters described by five dependent variables.

\subsection{Test Cases}
We consider six stochastic test cases, where we use the Smeaheia discretization described in section~\ref{sec:Smeaheia_model}. The test cases are presented in order of increasing complexity, but it should be noted that the simpler cases are really special cases of the more complex models where certain parameters are kept fixed to the mean values. This allows us both to study the fault upscaling in detail, and then to introduce the interaction between fault and formation uncertainty.

In Case I and II, the permeability in the Vette fault is assumed uncertain and follows the distributions depicted in Figure~\ref{fig:SGR_dist_2}. For practical sampling, lognormal distributions are fitted to the data for the different settings of $\perm[\text{c}]$. In Case I all other  parameters including the two-phase flow functions are assumed deterministic and set to their respective expected values (red curves in Figure~\ref{fig:param_param}), i.e., nominal reference values. In Case II we assume that each of the six layers of the reservoir model has lognormal homogeneous permeability, given by $k_i \sim \text{Lognormal}(\mu_{k_i}, \sigma_{k_i})$ for $i=1,\dots, 6$. The distribution parameters $\mu_{k_i}, \sigma_{k_i}$ are chosen so that the expectations of the layer permeabilities are, respectively, $1000, 50, 1000, 50, 850, 25$ mD with the standard deviations for simplicity set to 100 mD for all layers. In total, this model has seven independent random input parameters.

For the remaining Cases III-VI, we assume stochastic two-phase flow functions  (permeability, capillary pressure, saturation and wetting/non-wetting relative permeabilities) as represented by the copula model. Case III assumes all other parameters to be deterministic, while Case IV
is supplemented by an upscaled model for the uncertain permeability in the connection with Troll. This permeability is assumed to be independent of the Vette fault parameters, and computed using the same SGR-permeability mapping Eq.~\eqref{eq:SGR-perm}. The mean SGR values used to populate the geocellular model are given by the synthetic data in the depth range from 1500 to 1600 m shown in Figure~\ref{fig:SGR_data}, and the standard deviations are assumed to be 14 \% as for the Vette fault. Arithmetic averaging over the facies is performed to accommodate for flow across the fault. 

In Case V we consider the upscaled flow functions with the copula model in combination with the six reservoir layers following the distributions described in connection with Case II. Finally, in Case VI, we combine the random models from Case V (upscaled flow functions and reservoir layers) with uncertainty in the Troll connection permeability with the model described in Case IV.
The setups of the test cases are summarized in Table~\ref{tab:test_cases_summary}.

\begin{table}
    \centering
    \begin{adjustbox}{max width=\textwidth}
   \begin{tabular}{c|l | c}
     Case & Stochastic model & $n$ \\
     \hline
     I  & Fault permeability $K$ from lognormal data fit & 1\\
     \hline
     II  & 
     \begin{tabular}{@{}l@{}}
     Fault permeability $K$ as in Case I
     \\ Layer permeabilities $k_i$ ($i=1,\dots, 6$) lognormal with parameters s.t.\\
     $\E{k_i} = \{1000, 50, 1000, 50, 850, 25\}$ mD, $\text{Std}(k_i)=100$ mD.
     \end{tabular} & 7\\
     \hline
   III  & $K$, $\upscaledPc(\sdisc)$, $\upscaledsat(\sdisc)$, $K_{\text{r}}^{w}(\sdisc), K_{\text{r}}^{\text{nw}}(\sdisc)$  from copula + Eq.~\eqref{eq:K_Y1}-\eqref{eq:Knw_Y5} & 5 \\
   \hline
   IV  & 
   \begin{tabular}{@{}l@{}}
   $K$, $\upscaledPc(\sdisc)$, $\upscaledsat(\sdisc)$, $K_{\text{r}}^{w}(\sdisc)$, $K_{\text{r}}^{\text{nw}}(\sdisc)$ as in Case III \\
   $K_{\text{Troll}}$ Lognormal data fit
   \end{tabular} & 6\\
   \hline
   V  & 
   \begin{tabular}{@{}l@{}}
   Layer permeabilities as in Case II\\
   $K$, $\upscaledPc(\sdisc)$, $\upscaledsat(\sdisc)$, $K_{\text{r}}^{w}(\sdisc)$, $K_{\text{r}}^{\text{nw}}(\sdisc)$ as in Case III 
   \end{tabular} & 11
   \\
   \hline
   VI &
   \begin{tabular}{@{}l@{}}
   Layer permeabilities as in Case II\\
   $K$, $\upscaledPc(\sdisc)$, $\upscaledsat(\sdisc)$, $K_{\text{r}}^{w}(\sdisc)$, $K_{\text{r}}^{\text{nw}}(\sdisc)$ as in Case III \\
   $K_{\text{Troll}}$ as in Case IV
   \end{tabular}
   & 12
    \end{tabular}
    \end{adjustbox}
    \caption{Summary of the test cases. Only the stochastic parameters and functions are stated.  All other quantities are assumed fixed at their deterministic reference values. The number of stochastic dimensions denoted by $n$.}
    \label{tab:test_cases_summary}
\end{table}

\subsection{Monte Carlo Simulation}
\begin{figure}
    \centering  \subfigure[$k_{\textup{c}}=0.0001$ mD.]
{\includegraphics[width=0.99\textwidth]{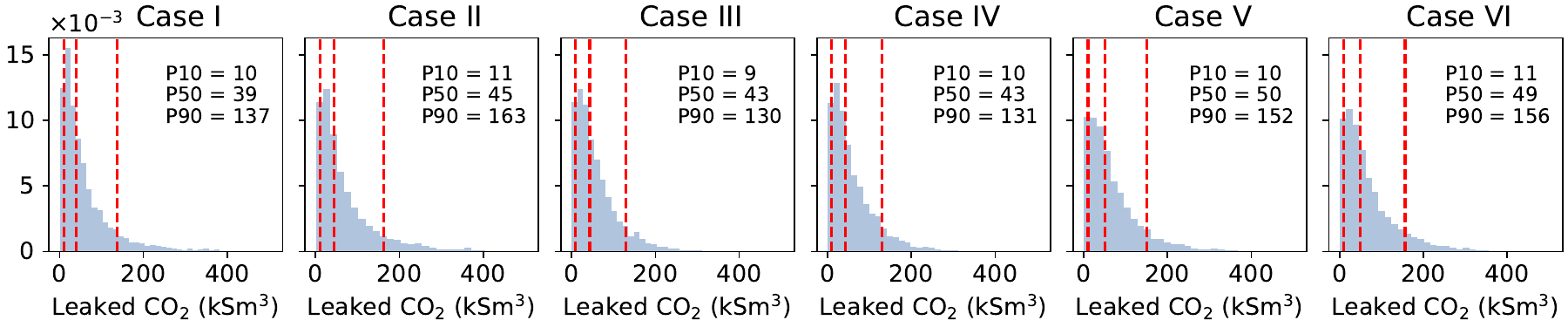}}
\subfigure[$k_{\textup{c}}=0.001$ mD.]
{\includegraphics[width=0.99\textwidth]{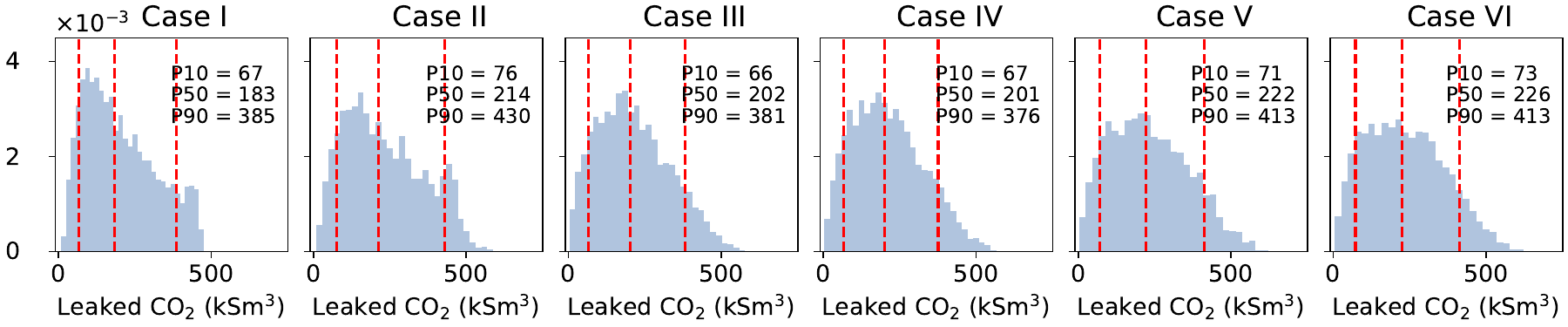}}
\subfigure[$k_{\textup{c}}=1$ mD.]
{\includegraphics[width=0.99\textwidth]{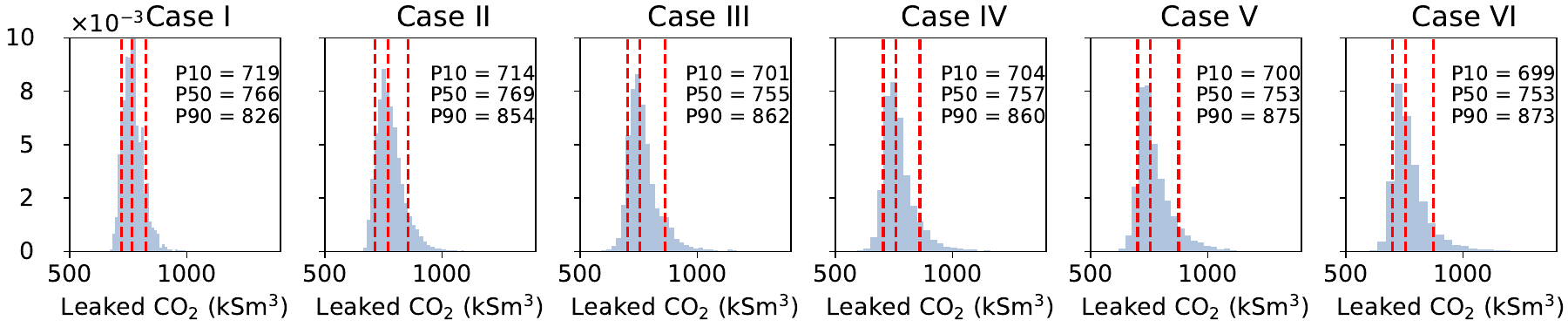}}
    \caption{Histograms of leaked CO$_2$ for Case I-VI. The three vertical dashed lines in each subfigure indicate the 10\textsuperscript{th} percentile (left), median (middle), and 90\textsuperscript{th} percentile (right).} 
    \label{fig:leakage_hist}
\end{figure}

The leakage distributions for the six test cases estimated from 5000 SMC samples and visualized as histograms are shown in 
Figure~\ref{fig:leakage_hist} for $\perm[\text{c}]=0.0001$ mD, $\perm[\text{c}]=0.001$ mD and $\perm[\text{c}]=1$ mD, respectively. As expected, the lower shale permeability leads to less amount of CO$_2$ leaked. 
The probability that there is no leakage at all is below 0.01 for the case $\perm[\text{c}]=0.0001$ mD, and identically zero for larger values of $\perm[\text{c}]$. For all clay permeability cases, the setup Case I differs from Cases II-VI in capturing the upper tail of the probability distributions. As a quantitative characterization of the leakage distribution, each figure includes the percentiles P10, P50 (median), and P90. More extreme percentiles (e.g., P99) are not included as they would not be accurately estimated with the current number of samples.
The variability of the distributions is larger between the results for different $\perm[c]$ than between the six cases of different stochastic models. This suggests that improved SGR models, or modeling $\perm[\text{c}]$ as a random variable, may yield more accurate predictions of the leakage distributions. The limited variability between the six cases of different stochastic models indicate that the simpler models are sufficient to obtain the statistics investigated for the current problem setup. However, a physically more refined model, or some other choice of statistics (e.g., rare events), may benefit from the more complex models involving the uncertainty in two-phase flow functions.

\subsection{Adaptive Stratified Sampling}
\begin{figure}
    \centering 
{\includegraphics[width=0.99\textwidth]{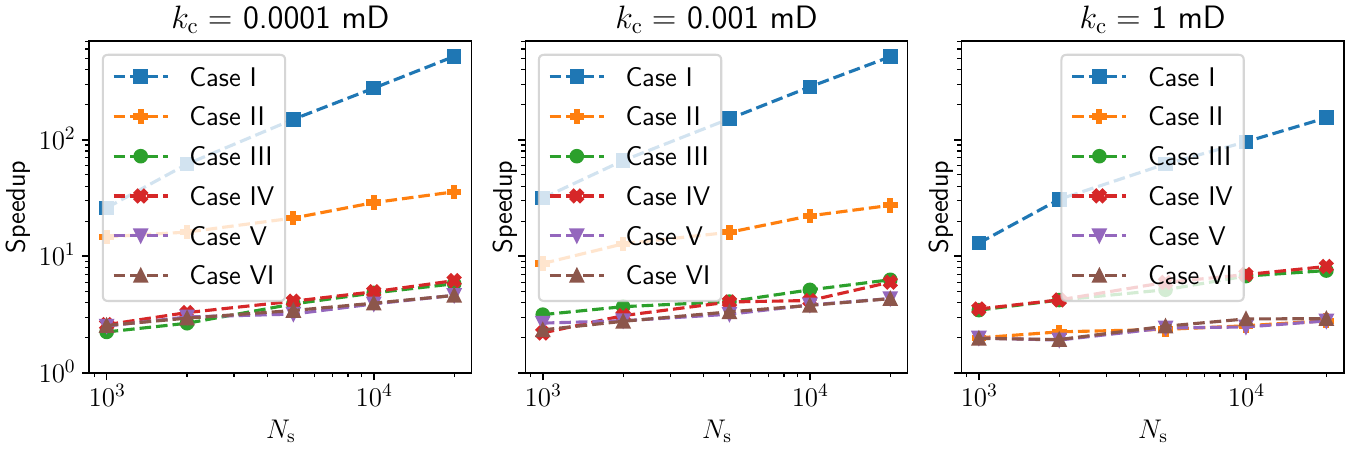}}
    \caption{Estimated speedup of ADSS with respect to SMC.} 
    \label{fig:speedups_Case_I_VI}
\end{figure}

Next we apply the ADSS method described in section~\ref{sec:adaptive_sampling} to the six test cases to compute the expected leakage using~\eqref{eq:ss_estimator:hybrid} with the hybrid sampling parameter $\alpha=0.5$, and allocation of 50 samples per iteration. These design parameter settings were found to be suitable for different test problems in~\cite{Pettersson_Krumscheid_22}, but have not been tailored to the test problems in the current paper.
The speedups of the ADSS method compared to SMC as given by~\eqref{eq:speedup} for all test cases and varying number of sampling budgets are shown in Figure~\ref{fig:speedups_Case_I_VI}.
For the simplest Case I, we observe two orders of magnitude speedup. For Case II with seven independent random variables, the speedup is one order of magnitude for the SGR models with $\perm[\text{c}]=0.0001, 0.001$ mD, but significantly lower for $\perm[\text{c}]=1$ mD. For all other cases, the speedups vary from 2-4 for the smallest sample sizes, to 2-8 for the largest sample sizes. The lower speedups are observed for the Cases III-VI where the relatively complex two-phase model~\eqref{eq:K_Y1}-\eqref{eq:Knw_Y5} is employed, with its sample path crossings of $\upscaledsat$ and $K_{\text{r}}^{w}$, and random variable dependence described by the Vine copula model as discussed in section~\ref{sec:copula_fit}. For $\perm[c]=1$ mD, we observe higher speedups for Cases III-IV than V-VI, which is to be expected given the difference in the number of random dimensions. Any localization of variance is more difficult to resolve for ADSS in high-dimensional spaces.
Taken together, the speedup results indicate that  the increased complexity of ADSS is motivated by the reduction in sample sizes possible while maintaining similar accuracy as SMC.

\subsection{Discussion} 

The current work presents a framework for uncertainty quantification for complex physical systems, e.g., subsurface \co storage, where the focus is on accurate and efficient representation of uncertainty from available data, and a subsequent uncertainty propagation through the physical model via adaptive sampling. If some input parameters to the physical model are considered uncertain but lacking in data, they should nevertheless be included as random variables, e.g., uniform random variables on an interval including possible physical variation. 
The proposed uncertainty modeling and propagation framework has been applied to a test case from the Smeaheia formation in the North Sea, but should be applicable to a wider class of problems.
Uncertainty in the physical model itself (in the current case the reservoir model) has not been considered in this work. Below we discuss separate aspects of the work presented.

A geocellular model with uncertain properties modelled as random variables has been upscaled, resulting in stochastic effective flow functions.
All upscaled flow functions except fault permeability are functions of a deterministic parameter and need to be modeled as random fields, in addition to the copula model used to represent the dependencies at some fixed value of the deterministic parameter. Karhunen-Loeve expansions, commonly used for this purpose, fail to capture physical features such as monotonic sample paths. Instead, we use data to fit an empirical correlation model that does not assume stochastic stationarity, i.e., the two-point correlation is allowed to change with the deterministic model parameter. All observed realizations of the resulting model are physically meaningful, although we cannot prove that the model excludes unphysical results.

Capturing the dependencies present in the multivariate distributions representing the upscaled model flow functions rely on existing methodology (and practically speaking software) for multidimensional copula fitting, effectively decomposing the problem into separate treatment of marginal distributions and dependence structure.
The considered copula approach offers flexibility, providing a good number of distinct pair copulas that can be fitted to high accuracy.
In particular, PDFs concentrated in the vicinity of the physical boundary of a parameter range (relative permeabilities and saturation are all confined to the unit interval) are well represented without any visible spurious boundary effects.
Indeed, the fitted copula models (for three different values of clay permeability) share the same tree structures, but the bivariate copula functions connecting the nodes differ somewhat, reflecting the differences in pairwise PDFs that have been observed. There are no intuitive and standardized measures for assessing proximity between multivariate distributions, but examining the 1D and 2D projections of the data and the model at least give an idea of the model fit, and they do agree well.

In addition to accurately representing empirical distributions of data, the upscaled model allows targeted model evaluation, i.e., directed sampling within a given region of random space. 
In standard Monte Carlo sampling, it is sufficient to be able to draw any random samples as long as they follow prescribed distributions honoring data. The adaptivity of the adaptive stratified sampling method employed here in addition requires that the stochastic input domain can be sampled from a given subdomain. Such subdomain sampling in a reduced and uniform space is made possible by the copula model and inverse Rosenblatt transformation. Many other groups of uncertainty propagation methods, including Quasi-Monte Carlo, Multi-Level Monte Carlo, and spectral expansions (generalized polynomial chaos) would either require or greatly benefit from these properties.

Adaptive stratified sampling requires significantly fewer samples than SMC to reach a given error threshold and is relatively insensitive to the number of stochastic dimensions, in particular if the number of effective dimensions is small (i.e., the total variance is dominated by a small subset of random dimensions). 
The design parameters set by the user, i.e., the number of samples allocated in every iteration and the proportion of samples optimally allocated (as opposed to the more robust proportional allocation), will affect the observed speedups.
In this work, a relatively conservative approach has been used, with as many as on average 50 samples per stratum allocated every iteration, and only a fraction $\alpha=0.5$ of samples optimally allocated. While any exact recommendation for these design parameters is hard to motivate, the chosen numbers appear suitable given the complex features of the Cases III-VI, which are illustrated in Figure~\ref{fig:flow_funs_vs_Us}.

Accelerated methods including adaptive stratified sampling provide improvement over SMC with respect to some user-defined metric, e.g., a mean value. SMC is hence more amenable to general analysis of quantities of interest, and was chosen to produce estimates of the distributions of \co leakage.
 The histograms of leaked CO$_2$ obtained by SMC indicate the possibility of long upper tails, i.e., small but nonzero probability for relatively large leakage. If these probabilities are of interest, a rare-event simulation framework may replace the current method for UP. The number of samples required is otherwise inversely proportional to the probability of the outcome of interest, so very unlikely scenarios are correspondingly difficult to correctly estimate.

Accurate uncertainty quantification relies on a similar level of accuracy of the physical model describing the complex system of interest. In the current work this includes both the fault model and the model for the surrounding reservoir.
The fault model considered in section \ref{sec:stochastic_faults} assumed a 1D description of the fine-scale geometry of the fault.
In reality, there can be significant horizontal variations both along and across the fault, thus in the general case it could be necessary to represent the fine-scale geometry as truly 3D. While we expect that the overall stochastic modeling framework to be still applicable, ingredients of the workflow will need to change: The upscaling must be based on numerical rather than analytical methods and since the upscaled stochastic ensembles may take different forms, the stochastic modeling approach could also need adaptations. 

The simulation framework successfully produced meaningful estimates of expected fault-related leakage for a prospective storage site in the Norwegian North Sea with relatively little computational overhead. It should thus be feasible to implement this approach in fault de-risking workflows used in an industrial context, i.e. it is not simply a theoretical exercise. One notable observation is that the leakage estimates obtained in the 18 simulated cases of the Smeaheia test case range from 80 to 1500 tons cumulative leaked \co over 60 years. This amounts to very tiny amounts of \co leakage that would be below the detection limits of modern subsea sensing equipment. We expect that these leakage estimates along the Vette fault test are indeed representative given that the underlying SGR map, reservoir model, and other parameters are based on actual data. However, more work needs to be done to ensure that the underlying assumptions in the dynamic simulations (i.e. NNC connections and grid resolution) are indeed valid. In addition, the results may be sensitive to the simplicity of the fault model (already discussed above) or require other stochastic simulation tools to capture the upper tails of the leakage distributions. In summary, the exact histograms and mean values derived in this study are for demonstration purposes and should not be taken out of context.



\section{Summary and Conclusions}
We have presented a numerical framework for upscaling of uncertain fault properties, uncertainty modeling including complexity reduction, and adaptive uncertainty propagation, and demonstrated numerical results for a real-world CO$_2$ storage test case.  
The combined physical and stochastic upscaling leads to a reduction in computational and random complexity, i.e., both decreased deterministic model evaluation time and smaller number of random dimensions. To exploit the stochastic model reduction, the reduced set of random variables are parameterized using an accurate copula model that captures both their dependence and their marginal distributions. The copula model is then combined with an empirical (data driven) model for the non-stationary correlation in deterministic parameter space, where standard methods such as Karhunen-Loeve expansion fail to produce physically relevant flow function models.

To reduce the number of physical model evaluations needed for some prescribed error tolerance of estimated mean values of leaked CO$_2$, ADSS has been used instead of direct SMC. Compared to SMC, for the sample size ranges that have been investigated in this paper, the computational budget for ADSS can be reduced with up to a factor of, respectively, 10 and 100 for the simpler stochastic cases, and 2--8 times for the copula models with nonlinear transformations. The speedups are sufficient to motivate the choice of using ADSS, rather than the more straightforward choice of SMC. In summary, the successful demonstration of the methodology are promising for future implementation of the framework in a real setting with more site-specific data available.

\section{Acknowledgments}
The work was funded by the Research Council of Norway through the projects Quantification of fault-related leakage risk (FRISK) under project number 294719, and Expansion of Resources for CO$_2$ Storage on the Horda Platform (ExpReCCS) under project number 336294.

\appendix
\section{Two-phase Upscaling}
\label{app:upscaling}
The two-phase upscaling applied in Section \ref{sec:twophase_upscaling} following \cite{rabinovich2016analytical} can be summarized as follows:
\begin{enumerate}
    \item Calculate an upscaled absolute permeability $K$.
    \item Pick a value for the capillary pressure, assign this to all fine-scale facies (this is the assumption of capillary equilibrium).
    \item Use an inverted fine-scale capillary pressure function to calculate a fine-scale saturation. If the fine-scale capillary pressure function is equal in all facies, the fine-scale saturation will be uniform.
    \item Calculate an upscaled saturation using arithmetic averaging.
    \item From the fine-scale saturation distribution, compute the upscaled effective phase-wise permeabilities, $\hat{K}^{\alpha}$, 
    as the harmonic average of the product of the fine-scale absolute and relative permeability, $k_{\text{r}}^{\alpha} k$. 
    \item Obtain the upscaled relative permeability as $K_{\text{r}}^{\alpha} = \hat{K}^{\alpha}/\hat{K}$ . 
\end{enumerate}
Points 2-6 are repeated over the full range of capillary pressure values;
we see that points 2 and 4 give an upscaled capillary pressure curve, while points 4 and 6 produce upscaled relative permeabilities.
In practice, point number 2 is a bit tricky since, first, the Brooks-Corey capillary pressure is not bounded at zero saturation, and second, the fine-scale capillary function is stochastic through its entry pressure, and varies between the fine-scale facies.
To achieve a robust implementation, the equilibrium values for capillary pressure are chosen as follows:
We sample a set of saturation values, denoted $\sdisc \in [1\times 10^{-6}, 1]$ where the cut-off is set slightly above 0 to avoid sampling an unbounded capillary pressure. The sampling is done in logarithmic space to achieve good coverage in the steep region of the Brooks-Corey curve. The equilibrium capillary pressure that enters point 2 is calculated from these values of $\sdisc$, using the minimal entry pressure of the facies in the current realization, thereby ensuring that the equilibrium pressure can be reached in all facies.

\section{Copula Modeling}
\subsection{Example: Conditioning on Bivariate Copulas}
\label{sec:cop-fact-example}
As an example of conditioning on bivariate copulas,  borrowed from~\cite{Czado_Nagler_22}, the trivariate joint PDF $f_{Y}(y_1, y_2, y_3)$ can be represented by

\begin{multline}
f_{Y}(y_1, y_2, y_3) =c_{1,3 | 2}(F_{1 | 2}(y_1 | y_2), F_{3 | 2}(y_3 | y_2) | y_2) c_{23}(F_2(y_2), F_{3}(y_3))\\        \times c_{12}(F_1(y_1), F_{2}(y_2)) f_{1}(y_1) f_{2}(y_2) f_{3}(y_3),
\end{multline}
but we also have the equally valid factorizations
\begin{align*}
f_{Y}(y_1, y_2, y_3) &=c_{2,3 | 1}(F_{2 | 1}(y_2 | y_1), F_{3 | 1}(y_3 | y_1) | y_1) c_{13}(F_1(y_1), F_{3}(y_3))\\        & \times c_{12}(F_1(y_1), F_{2}(y_2)) f_{1}(y_1) f_{2}(y_2) f_{3}(y_3), \\
f_{Y}(y_1, y_2, y_3) &= c_{1,2 | 3}(F_{1 | 3}(y_1 | y_3), F_{2 | 3}(y_2 | y_3) | y_3) c_{13}(F_1(y_1), F_{3}(y_3))\\        & \times c_{23}(F_2(y_2), F_{3}(y_3)) f_{1}(y_1) f_{2}(y_2) f_{3}(y_3),
\end{align*}
so that the representation is not unique.

\subsection{Copula Tree Fit}
\label{sec:cop-tree-fit}
\begin{figure}[H]

\subfigure[1 mD.] 
{ 
\begin{tikzpicture}[thick,scale=0.6, every node/.style={transform shape}]
    \Text[x=-1, y=3]{Tree 1}
    \Vertex[color=white, size=0.8, fontsize=\normalsize, x=-1, y=2, label=2]{1}
    \Vertex[color=white, size=0.8, fontsize=\normalsize, x=2, y=2, label=3]{2}
    \Vertex[color=white, size=0.8, fontsize=\normalsize, x=5, y=2, label=4]{3}
    \Vertex[color=white, size=0.8, fontsize=\normalsize, x=8, y=2, label=5]{4}
    \Vertex[color=white, size=0.8, fontsize=\normalsize, x=11, y=2, label=1]{5}
    \Edge[lw=1, label={2,3}, fontsize=\normalsize](1)(2)
    \Edge[lw=1, label={3,4}, fontsize=\normalsize](2)(3)
    \Edge[lw=1, label={4,5}, fontsize=\normalsize](3)(4)
    \Edge[lw=1, label={5,1}, fontsize=\normalsize](4)(5)
    \Text[x=0.5, y=1.5]{BB1 90$^{\circ}$}
    \Text[x=3.5, y=1.5]{Gumbel}
    \Text[x=6.5, y=1.5]{TLL}
    \Text[x=9.5, y=1.5]{TLL}
    
    \Text[x=-1, y=0]{Tree 2}
    \Vertex[shape=ellipse, color=white, size=0.9, fontsize=\normalsize, x=-1, y=-1, label={2,3}]{onetwo}
    \Vertex[color=white, size=0.9, fontsize=\normalsize, x=2, y=-1, label={3,4}]{twothree}
    \Vertex[color=white, size=0.9, fontsize=\normalsize, x=5, y=-1, label={4,5}]{threefour}
    \Vertex[color=white, size=0.9, fontsize=\normalsize, x=8, y=-1, label={5,1}]{fourfive}
    
    \Edge[lw=1, label={$2,4 | 3$}, fontsize=\normalsize](onetwo)(twothree)
    \Edge[lw=1, label={$3,5  |  4$}, fontsize=\normalsize](twothree)(threefour)
    \Edge[lw=1, label={$1,4 | 5$}, fontsize=\normalsize](threefour)(fourfive)
    
    \Text[x=0.5, y=-1.5]{TLL}
    \Text[x=3.5, y=-1.5]{TLL}
    \Text[x=6.5, y=-1.5]{TLL}

    \Text[x=-1, y=-3]{Tree 3}
    \Vertex[shape=ellipse, color=white, size=0.8, fontsize=\normalsize, x=-0.5, y=-4, label={$2,4\ | \ 3$}, style={minimum width=1.8cm}]{onethreetwo}
    \Vertex[shape=ellipse, color=white, size=0.8, fontsize=\normalsize, x=4, y=-4, label={$3,5 \ | \ 4$}, style={minimum width=1.8cm}]{twofourthree}
    \Vertex[shape=ellipse, color=white, size=0.8, fontsize=\normalsize, x=8.5, y=-4, label={$1,4 \ | \ 5$}, style={minimum width=1.8cm}]{threefivefour}
    
    \Edge[lw=1, label={$2,5 | 3,4$}, fontsize=\normalsize](onethreetwo)(twofourthree)
    \Edge[lw=1, label={$1,3 | 4,5$}, fontsize=\normalsize](twofourthree)(threefivefour)
    \Text[x=1.75, y=-4.5]{TLL}
    \Text[x=6.25, y=-4.5]{Gaussian}

    \Text[x=-1, y=-6]{Tree 4}
    \Vertex[shape=ellipse, color=white, size=0.8, fontsize=\normalsize, x=-0.5, y=-7, label={$2,5\ | \ 3,4$}, style={minimum width=2.2cm}]{onefourtwothree}
    \Vertex[shape=ellipse, color=white, size=0.8, fontsize=\normalsize, x=5, y=-7, label={$1,3 \ | \ 4,5$}, style={minimum width=2.2cm}]{twofivethreefour}
    
    \Edge[lw=1, label={$1,2 | 3,4,5$}, fontsize=\normalsize](onefourtwothree)(twofivethreefour)
    \Text[x=2.25, y=-7.5]{Student-t}
    
\end{tikzpicture}
}
~
\subfigure[$k_{\text{c}} =$0.001 mD.] 
{ 
\begin{tikzpicture}[thick,scale=0.6, every node/.style={transform shape}]
    \Text[x=-1, y=3]{Tree 1}
    \Vertex[color=white, size=0.8, fontsize=\normalsize, x=-1, y=2, label=2]{1}
    \Vertex[color=white, size=0.8, fontsize=\normalsize, x=2, y=2, label=3]{2}
    \Vertex[color=white, size=0.8, fontsize=\normalsize, x=5, y=2, label=4]{3}
    \Vertex[color=white, size=0.8, fontsize=\normalsize, x=8, y=2, label=5]{4}
    \Vertex[color=white, size=0.8, fontsize=\normalsize, x=11, y=2, label=1]{5}
    \Edge[lw=1, label={2,3}, fontsize=\normalsize](1)(2)
    \Edge[lw=1, label={3,4}, fontsize=\normalsize](2)(3)
    \Edge[lw=1, label={4,5}, fontsize=\normalsize](3)(4)
    \Edge[lw=1, label={5,1}, fontsize=\normalsize](4)(5)
    \Text[x=0.5, y=1.5]{TLL}
    \Text[x=3.5, y=1.5]{TLL}
    \Text[x=6.5, y=1.5]{Gaussian}
    \Text[x=9.5, y=1.5]{TLL}
    
    \Text[x=-1, y=0]{Tree 2}
    \Vertex[shape=ellipse, color=white, size=0.9, fontsize=\normalsize, x=-1, y=-1, label={2,3}]{onetwo}
    \Vertex[color=white, size=0.9, fontsize=\normalsize, x=2, y=-1, label={3,4}]{twothree}
    \Vertex[color=white, size=0.9, fontsize=\normalsize, x=5, y=-1, label={4,5}]{threefour}
    \Vertex[color=white, size=0.9, fontsize=\normalsize, x=8, y=-1, label={5,1}]{fourfive}
    
    \Edge[lw=1, label={$2,4 | 3$}, fontsize=\normalsize](onetwo)(twothree)
    \Edge[lw=1, label={$3,5  |  4$}, fontsize=\normalsize](twothree)(threefour)
    \Edge[lw=1, label={$1,4 | 5$}, fontsize=\normalsize](threefour)(fourfive)
    
    \Text[x=0.5, y=-1.5]{TLL}
    \Text[x=3.5, y=-1.5]{TLL}
    \Text[x=6.5, y=-1.5]{TLL}

    \Text[x=-1, y=-3]{Tree 3}
    \Vertex[shape=ellipse, color=white, size=0.8, fontsize=\normalsize, x=-0.5, y=-4, label={$2,4\ | \ 3$}, style={minimum width=1.8cm}]{onethreetwo}
    \Vertex[shape=ellipse, color=white, size=0.8, fontsize=\normalsize, x=4, y=-4, label={$3,5 \ | \ 4$}, style={minimum width=1.8cm}]{twofourthree}
    \Vertex[shape=ellipse, color=white, size=0.8, fontsize=\normalsize, x=8.5, y=-4, label={$1,4 \ | \ 5$}, style={minimum width=1.8cm}]{threefivefour}
    
    \Edge[lw=1, label={$2,5 | 3,4$}, fontsize=\normalsize](onethreetwo)(twofourthree)
    \Edge[lw=1, label={$1,3 | 4,5$}, fontsize=\normalsize](twofourthree)(threefivefour)
    \Text[x=1.75, y=-4.5]{Gaussian}
    \Text[x=6.25, y=-4.5]{Gaussian}

    \Text[x=-1, y=-6]{Tree 4}
    \Vertex[shape=ellipse, color=white, size=0.8, fontsize=\normalsize, x=-0.5, y=-7, label={$2,5\ | \ 3,4$}, style={minimum width=2.2cm}]{onefourtwothree}
    \Vertex[shape=ellipse, color=white, size=0.8, fontsize=\normalsize, x=5, y=-7, label={$1,3 \ | \ 4,5$}, style={minimum width=2.2cm}]{twofivethreefour}
    
    \Edge[lw=1, label={$1,2 | 3,4,5$}, fontsize=\normalsize](onefourtwothree)(twofivethreefour)
    \Text[x=2.25, y=-7.5]{Student-t}
    
\end{tikzpicture}
}
\caption{Tree structure of the fitted copulas. This tree structure is known as a D-vine copula.}
\label{fig:vine-trees}
\end{figure}
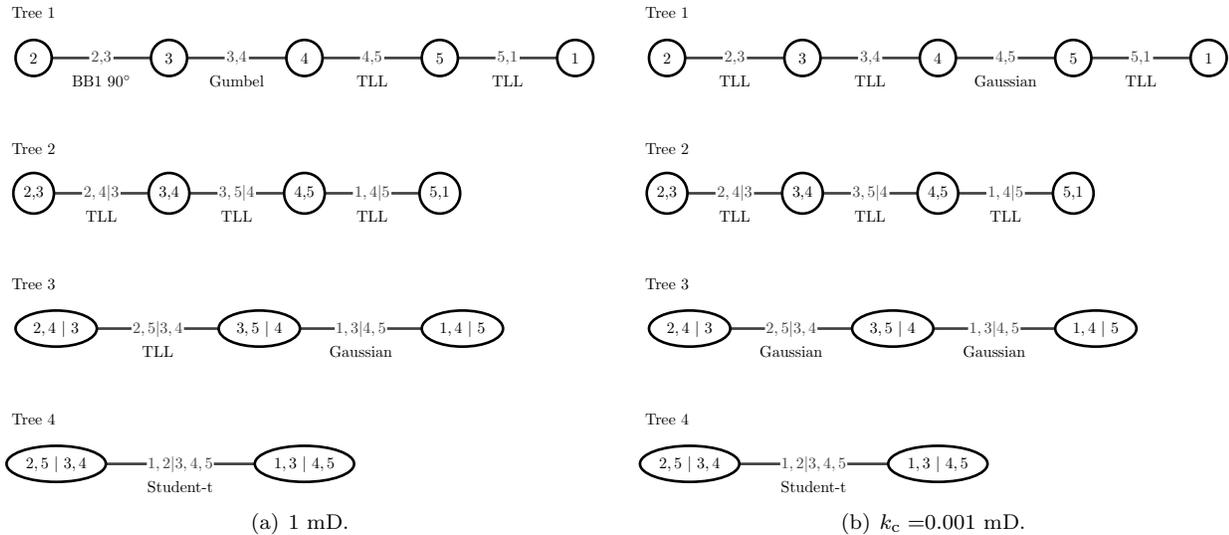

The tree structures defining the fitted copulas are shown in Fig.~\ref{fig:vine-trees}. The nodes in the trees denote random variables, where the random variables are represented by their indices (e.g., $1,2$ is short for $Y_1,Y_2$). The edge connecting two adjacent nodes denotes the joint distribution conditioned on the subset of repeated variables of the two nodes. The edges of Tree $t$ are the nodes of Tree $t+1$, so in this sense the structure of all subsequent trees can be inferred from Tree 1. The choice of copula function, among a suitable candidate set of copula families, for any pairs of variables is inferred by data.
The PDF of the five-dimensional distribution is the product of all (ten) bivariate distributions represented by the edges of the trees in Fig.~\ref{fig:vine-trees} and the (five) marginal distributions. The bivariate copulas that appear in Fig.~\ref{fig:vine-trees} are the nonparametric transformation local likelihood (TLL) copula~\cite{Geenens_etal_17}, and common parametric function copulas (BB1~\cite{Joe_97}, Gumbel, Student-t, and Gaussian) that may in some cases be rotated by 90$^{\circ}$, 180$^{\circ}$, or 270$^{\circ}$ for best fit with data. The copula model of the case $\perm[c]=0.0001$ mD has the same tree structure as shown in Fig.~\ref{fig:vine-trees}, but with bivariate copulas (from left to right) as follows. Tree 1: TLL, BB7 180, Gaussian, TLL; Tree 2: TLL, TLL, TLL; Tree 3: TLL, BB8; Tree 4: BB7 270$^{\circ}$.

\bibliographystyle{plain}
\bibliography{bibliography.bib}

\begin{thebibliography}{10}

\bibitem{Aas_Berg_09}
K.~Aas and D.~Berg.
\newblock Models for construction of multivariate dependence – a comparison
  study.
\newblock {\em Eur. J. Financ.}, 15(7-8):639--659, 2009.

\bibitem{opmFlowManual}
J.~Alvestad, D.~Baxendale, K.~Bao, M.~Blatt, J.~Hove, A.~Lauser, C.~Goncalves
  Machado, D.~Landa-Marb\'{a}n, A.~M. Kvarving, A.~Rasmussen, A.~B. Rustad,
  T.~H. Sandve, B.~Skaflestad, T.~Skille, and P.~J. Verveer.
\newblock {\em {OPM Flow Reference Manual (2021-10)}}.
\newblock {Open Porous Media Initiative}.
\newblock \url{https://opm-project.org/?page\_id=955}.

\bibitem{Aziz_Settari_79}
K.~Aziz and A.~Settari.
\newblock {\em Petroleum Reservoir Simulation}.
\newblock Springer Netherlands, 1979.

\bibitem{Babuska_etal_07}
I.~Babu\v{s}ka, F.~Nobile, and R.~Tempone.
\newblock A stochastic collocation method for elliptic partial differential
  equations with random input data.
\newblock {\em {SIAM} J. Numer. Anal.}, 45(3):1005--1034, 2007.

\bibitem{Bedford_Cooke_02}
T.~Bedford and R.~M. Cooke.
\newblock Vines--a new graphical model for dependent random variables.
\newblock {\em Ann. Stat.}, 30(4):1031--1068, 2002.

\bibitem{Berge_etal_22}
Runar~L Berge, Sarah~E Gasda, Eirik Keilegavlen, and Tor~Harald Sandve.
\newblock Impact of deformation bands on fault-related fluid flow in
  field-scale simulations.
\newblock {\em Int. J. Greenh. Gas Con.}, 119:103729, 2022.

\bibitem{Bhatti_Do_19}
M.~I. Bhatti and H.~Q. Do.
\newblock Recent development in copula and its applications to the energy,
  forestry and environmental sciences.
\newblock {\em Int. J. Hydrogen Energ.}, 44(36):19453--19473, 2019.

\bibitem{Bjornaraa_etal_21}
T.~I. Bj{\o}rnar\aa, E.~M. Haines, and E.~Skurtveit.
\newblock Upscaled geocellular flow model of potential across-and along-fault
  leakage using shale gouge ratio.
\newblock In {\em proceedings of TCCS, Trondheim}, 2021.

\bibitem{Bjornaraa_etal_22}
T.~I. Bj{\o}rnar\aa, E.~Skurtveit, E.~Michie, and S.~A. Smith.
\newblock Overburden fluid migration along the {Vette Fault Zone, North Sea},
  using different fault permeability models.
\newblock 2022(1):1--5, 2022.

\bibitem{Braathen_etal_09}
A.~Braathen, J.~Tveranger, H.~Fossen, T.~Skar, N.~Cardozo, S.E. Semshaug,
  E.~Bastesen, and E.~Sverdrup.
\newblock Fault facies and its application to sandstone reservoirs.
\newblock {\em AAPG bulletin}, 93(7):891--917, 2009.

\bibitem{Brechmann_Schepsmeier_13}
E.~C. Brechmann and U.~Schepsmeier.
\newblock Modeling dependence with {C}- and {D}-vine copulas: The {R} package
  {CDVine}.
\newblock {\em J. Stat. Softw.}, 52(3):1--27, 2013.

\bibitem{brooks1965hydraulic}
Royal~Harvard Brooks.
\newblock {\em Hydraulic properties of porous media}.
\newblock Colorado State University, 1965.

\bibitem{Childs_etal_90}
C.~Childs, J.~J. Walsh, and J.~Watterson.
\newblock A method for estimation of the density of fault displacements below
  the limits of seismic resolution in reservoir formations.
\newblock In A.~T. Buller, E.~Berg, O.~Hjelmeland, J.~Kleppe, O.~Tors{\ae}ter,
  and J.~O. Aasen, editors, {\em North {S}ea Oil and Gas Reservoirs---II},
  pages 309--318. Springer Netherlands, 1990.

\bibitem{Czado_Nagler_22}
C.~Czado and T.~Nagler.
\newblock Vine copula based modeling.
\newblock {\em Annu. Rev. Stat. Appl.}, 9:453--477, 2022.

\bibitem{Elsheikh_etal_14}
A.~H. Elsheikh, S.~Oladyshkin, W.~Nowak, and M.~Christie.
\newblock Estimating the probability of {CO}$_2$ leakage using rare event
  simulation.
\newblock In {\em {ECMOR XIV}-14th European conference on the mathematics of
  oil recovery}, volume 2014, pages 1--9. European Association of Geoscientists
  \& Engineers, 2014.

\bibitem{Etore_etal_11}
P.~{\'E}tor{\'e}, G.~Fort, B.~Jourdain, and E.~Moulines.
\newblock On adaptive stratification.
\newblock {\em Ann. Oper. Res.}, 89(1):127--154, 2011.

\bibitem{UNFCCC}
{European Union (Convention)}.
\newblock {UNFCCC} greenhouse gas inventory data, 2022.
\newblock Available at \url{https://di.unfccc.int/flex\_annex1}.

\bibitem{Fisher_Jolley_07}
Q.~J. Fisher and S.~J. Jolley.
\newblock Treatment of faults in production simulation models.
\newblock {\em Geol. Soc. Spec. Publ.}, 292(1):219--233, 2007.

\bibitem{flemisch2011dumux}
B.~Flemisch, M.~Darcis, K.~Erbertseder, B.~Faigle, A.~Lauser, K.~Mosthaf,
  S.~M{\"u}thing, P.~Nuske, A.~Tatomir, M.~Wolff, et~al.
\newblock Dumux: Dune for multi-$\{$phase, component, scale, physics,…$\}$
  flow and transport in porous media.
\newblock {\em Adv. Water Res.}, 34(9):1102--1112, 2011.

\bibitem{Freeman_etal_08}
S.~R. Freeman, S.~D. Harris, and R.~J. Knipe.
\newblock Fault seal mapping – incorporating geometric and property
  uncertainty.
\newblock {\em Geol. Soc. Spec. Publ.}, 309(1):5--38, 2008.

\bibitem{Geenens_etal_17}
G.~Geenens, A.~Charpentier, and D.~Paindaveine.
\newblock Probit transformation for nonparametric kernel estimation of the
  copula density.
\newblock {\em Bernoulli}, 23(3):1848--1873, 2017.

\bibitem{Giles_2015}
M.~B. Giles.
\newblock Multilevel {M}onte {C}arlo methods.
\newblock {\em Acta Numerica}, 24:259--328, 2015.

\bibitem{gross2021geosx}
H.~Gross and A.~Mazuyer.
\newblock Geosx: A multiphysics, multilevel simulator designed for exascale
  computing.
\newblock In {\em SPE Reservoir Simulation Conference}. OnePetro, 2021.

\bibitem{Halland_etal_11}
E.~K. Halland, W.~T. Johansen, and F.~Riis.
\newblock {CO}$_2$ storage atlas {Norwegian North Sea}.
\newblock {\em Norwegian Petroleum Directorate, PO Box}, 600, 2011.

\bibitem{Joe_96}
H.~Joe.
\newblock Families of $m$-variate distributions with given margins and
  $m(m-1)/2$ bivariate dependence parameters.
\newblock In L.~R{\"u}schendorf, B.~Schweizer, and M.~D. Taylor, editors, {\em
  Distributions with Fixed Marginals and Related Topics}, volume~28, pages
  120--141, Hayward, CA, 1996. Institute of Mathematical Statistics.

\bibitem{Joe_97}
H.~Joe.
\newblock {\em Multivariate models and multivariate dependence concepts}.
\newblock CRC press, 1997.

\bibitem{Jolley_etal_07}
S.~J. Jolley, H.~Dijk, J.~H. Lamens, Q.~J. Fisher, T.~Manzocchi, H.~Eikmans,
  and Y.~Huang.
\newblock Faulting and fault sealing in production simulation models: Brent
  province, northern {North Sea}.
\newblock {\em Petrol. Geosci.}, 13(4):321--340, 2007.

\bibitem{kang2014analytical}
M.~Kang, J.~M. Nordbotten, F.~Doster, and M.~A. Celia.
\newblock Analytical solutions for two-phase subsurface flow to a leaky fault
  considering vertical flow effects and fault properties.
\newblock {\em Water Resour. Res.}, 50(4):3536--3552, 2014.

\bibitem{Kloek_vanDijk_78}
T.~Kloek and H.~K. {van Dijk}.
\newblock Bayesian estimates of equation system parameters: An application of
  integration by {Monte Carlo}.
\newblock {\em Econometrica}, 46(1):1--19, 1978.

\bibitem{Kolyukhin_Tveranger_15}
D.~Kolyukhin and J.~Tveranger.
\newblock Statistical modelling of fault core and deformation band structure in
  fault damage zones.
\newblock 06 2015.

\bibitem{Krevor_etal_23}
S.~Krevor, H.~de~Coninck, S.E. Gasda, and others.
\newblock Subsurface carbon dioxide and hydrogen storage for a sustainable
  energy future.
\newblock {\em Nat. Rev. Earth Environ.}, 2023.

\bibitem{Krumscheid_etal_package_23}
S.~Krumscheid, F.~Yuan, and P.~Pettersson.
\newblock Package '{ADSS}'.
\newblock 2023.

\bibitem{Krumscheid_Pettersson_23}
Sebastian Krumscheid and Per Pettersson.
\newblock Sequential estimation using hierarchically stratified domains with
  {L}atin hypercube sampling, 2023.

\bibitem{LEcuyer_Lemieux_02}
P.~L'Ecuyer and C.~Lemieux.
\newblock Recent advances in randomized {Quasi-Monte Carlo} methods.
\newblock In M.~Dror, P.~L'Ecuyer, and F.~Szidarovszky, editors, {\em Modeling
  Uncertainty: An Examination of Stochastic Theory, Methods, and Applications},
  International Series in Operations Research \& Management Science, pages
  419--474. Springer, New York, NY, 2002.

\bibitem{leverett1941capillary}
MoC Leverett.
\newblock Capillary behavior in porous solids.
\newblock {\em Transactions of the AIME}, 142(01):152--169, 1941.

\bibitem{lie2019introduction}
K.-A. Lie.
\newblock {\em An introduction to reservoir simulation using {MATLAB/GNU
  O}ctave: User guide for the {MATLAB} Reservoir Simulation Toolbox ({MRST})}.
\newblock Cambridge University Press, 2019.

\bibitem{Manzocchi_etal_08}
T.~Manzocchi, A.~E. Heath, B.~Palananthakumar, C.~Childs, and J.~J. Walsh.
\newblock Faults in conventional flow simulation models: a consideration of
  representational assumptions and geological uncertainties.
\newblock {\em Petrol. Geosci.}, 14(1):91--110, 2008.

\bibitem{McKay_etal_79}
M.~D. Mc{K}ay, R.~J. Beckman, and W.~J. Conover.
\newblock A comparison of three methods for selecting values of input variables
  in the analysis of output from a computer code.
\newblock {\em Technometrics}, 21(2):239--245, 1979.

\bibitem{Michie_etal_21}
E.~A.~H. Michie, M.~J. Mulrooney, and A.~Braathen.
\newblock Fault interpretation uncertainties using seismic data, and the
  effects on fault seal analysis: a case study from the {H}orda {P}latform,
  with implications for {CO}$_{2}$ storage.
\newblock {\em Solid Earth}, 12(6):1259--1286, 2021.

\bibitem{Nagler_etal_22package}
T.~Nagler, U.~Schepsmeier, J.~Stoeber, E.~C. Brechmann, B.~Graeler, T.~Erhardt,
  C.~Almeida, A.~Min, C.~Czado, M.~Hofmann, et~al.
\newblock Package ‘vinecopula’.
\newblock 2022.

\bibitem{Pei_etal_15}
Y.~Pei, D.~A. Paton, R.~J. Knipe, and K.~Wu.
\newblock A review of fault sealing behaviour and its evaluation in
  siliciclastic rocks.
\newblock {\em Earth-Sci. Rev.}, 150:121--138, 2015.

\bibitem{Perrin_13}
G.~Perrin, C.~Soize, D.~Duhamel, and C.~Funfschilling.
\newblock Karhunen–{L}oeve expansion revisited for vector-valued random
  fields: Scaling, errors and optimal basis.
\newblock {\em J. Comput. Phys.}, 242:607--622, 2013.

\bibitem{Pettersson_Krumscheid_22}
P.~Pettersson and S.~Krumscheid.
\newblock Adaptive stratified sampling for nonsmooth problems.
\newblock {\em Int. J. Uncertain. Quan.}, 12(6):71--99, 2022.

\bibitem{Pettersson_etal_22}
P.~Pettersson, S.~Tveit, and S.~E. Gasda.
\newblock Dynamic estimates of extreme-case {CO}$_2$ storage capacity for
  basin-scale heterogeneous systems under geological uncertainty.
\newblock {\em Int. J. Greenh. Gas Con.}, 116:103613, 2022.

\bibitem{pruess1991tough2}
K.~Pruess.
\newblock {TOUGH2} - a general-purpose numerical simulator for multiphase fluid
  and heat flow.
\newblock 1991.

\bibitem{rabinovich2016analytical}
A.~Rabinovich, B.~Li, and L.~J. Durlofsky.
\newblock Analytical approximations for effective relative permeability in the
  capillary limit.
\newblock {\em Water Resour. Res.}, 52(10):7645--7667, 2016.

\bibitem{Rahman_etal_21}
M.~J. Rahman, J.~C. Choi, M.~Fawad, and N.~H. Mondol.
\newblock Probabilistic analysis of {V}ette fault stability in potential
  {CO}$_2$ storage site {S}meaheia, offshore {N}orway.
\newblock {\em Int. J. Greenh. Gas Con.}, 108:103315, 2021.

\bibitem{rasmussen2021open}
A.~Fl{\o} Rasmussen, T.~H. Sandve, K.~Bao, A.~Lauser, J.~Hove, B.~Skaflestad,
  R.~Kl{\"o}fkorn, M.~Blatt, A.~B. Rustad, O.~S{\ae}vareid, et~al.
\newblock The open porous media flow reservoir simulator.
\newblock {\em Comput. Math. Appl.}, 81:159--185, 2021.

\bibitem{Ringrose_Meckel_19}
P.~S. Ringrose and T.~A Meckel.
\newblock Maturing global {CO}$_2$ storage resources on offshore continental
  margins to achieve {2DS} emissions reductions.
\newblock {\em Sci Rep}, 9(1):17944, 2019.

\bibitem{Rosenblatt_52}
M.~Rosenblatt.
\newblock Remarks on a multivariate transformation.
\newblock {\em Ann. Math. Stat.}, 23(3):470--472, 1952.

\bibitem{sandve2021convective}
T.~H. Sandve, S.~E. Gasda, A.~Rasmussen, and A.~B. Rustad.
\newblock Convective dissolution in field scale {CO}$_2$ storage simulations
  using the {OPM} flow simulator.
\newblock In {\em {TCCS}--11. {CO}$_2$ Capture, Transport and Storage.
  Trondheim 22nd--23rd June 2021 Short Papers from the 11th International
  {T}rondheim {CCS} Conference}. SINTEF Academic Press, 2021.

\bibitem{sandve2022simulators}
T.~H. Sandve, A.~B. Rustad, A.~Thune, B.~Nazarian, S.~Gasda, and A.~F.
  Rasmussen.
\newblock Simulators for the gigaton storage challenge. {A} benchmark study on
  the regional {S}meaheia model.
\newblock In {\em {EAGE} GeoTech 2022 Sixth {EAGE} Workshop on {CO}$_2$
  Geological Storage}, volume 2022, pages 1--5. European Association of
  Geoscientists \& Engineers, 2022.

\bibitem{Schueller_etal_13}
S.~Schueller, A.~Braathen, H.~Fossen, and J.~Tveranger.
\newblock Spatial distribution of deformation bands in damage zones of
  extensional faults in porous sandstones: Statistical analysis of field data.
\newblock {\em J. Struct. Geol.}, 52:148--162, 2013.

\bibitem{Shields_16}
M.~D. Shields.
\newblock Refined {L}atinized stratified sampling: A robust sequential sample
  size extension methodology for high-dimensional {L}atin hypercube and
  stratified designs.
\newblock {\em Int. J. Uncertain. Quan.}, 6(1):79--97, 2016.

\bibitem{Shields_etal_15}
M.~D. Shields, K.~Teferra, A.~Hapij, and R.~P. Daddazio.
\newblock Refined stratified sampling for efficient {M}onte {C}arlo based
  uncertainty quantification.
\newblock {\em Reliab. Eng. Syst. Safe.}, 142:310--325, 2015.

\bibitem{Sklar_59}
M.~Sklar.
\newblock Fonctions de r\'{e}partition \`{a} $n$ dimensions et leurs marges.
\newblock {\em Publ. inst. statist. univ. Paris}, 8:229--231, 1959.

\bibitem{Sperrevik_etal_02}
S.~Sperrevik, P.A. Gillespie, Q.J. Fisher, T.~Halvorsen, and R.~Knipe.
\newblock Empirical estimation of fault rock properties.
\newblock In A.~G. Koestler and R.~Hunsdale, editors, {\em Hydrocarbon Seal
  Quantification}, volume~11 of {\em Norwegian Petroleum Society (NPF), Special
  Publications}, page 109–125. 2002.

\bibitem{Wu_etal_21}
L.~Wu, R.~Thorsen, S.~Ottesen, R.~Meneguolo, K.~Hartvedt, P.~Ringrose, and
  B.~Nazarian.
\newblock Significance of fault seal in assessing {CO}$_2$ storage capacity and
  containment risks - an example from the {H}orda {P}latform, northern {N}orth
  {S}ea.
\newblock {\em Petrol. Geosci.}, 27(3):petgeo2020--102, 2021.

\bibitem{Xiu_Karniadakis_02}
D.~Xiu and G.~Em. Karniadakis.
\newblock The {W}iener--{A}skey polynomial chaos for stochastic differential
  equations.
\newblock {\em {SIAM} J. Sci. Comput.}, 24(2):619--644, 2002.

\bibitem{Yielding_etal_97}
G.~Yielding, B.~Freeman, and D.~T. Needham.
\newblock Quantitative fault seal prediction.
\newblock {\em AAPG Bulletin}, 81(6):897--917, 1997.

\end{thebibliography}

\end{document}